%% file: notesonemptylatticepolytopes
\input amstex
\input generic_macros

\def\flo #1{\left\lfloor {#1} \right\rfloor}

\def\ceil #1{\left\lceil {#1} \right\rceil}

\let\GL=\gl%

 %

 %

\def\Leg#1,#2.{\(\frac{#1}{{#2}}\)}
\def\Ip #1,#2.{\langle\!\langle #1,#2\rangle\! \rangle}
\def\Ipd #1,#2.{\Ip #1,#2._d}

\def\1{{\pmb 1}}

\def\cvx{\text{cvx\,}}

\def\brm #1,#2.{{}_{#1\backslash}#2}
\def\ah{asymptotically hollow{}}

\let\cont=\Ct


\NoBlackBoxes

\def\notdivide
{\,|{\hskip-.9ex}\backslash }

\def\SSS{5}
\def\aa{\alpha}

\def\cont{\text{cont\,}}
\def\Fr #1{\brcs{#1}} 

\def\agl #1{\text{AGL($#1$)}}
\def\wid{\text{wid\,}}

\def\oneone{1.1}
\def\onetwo{1.2}
\def\onethr{1.3}

\def\throne{3.1}
\def\thrtwo{3.2}
\def\thrthr{3.3}

\def\fivfiv{2.5}
\def\fivsix{2.6}
\def\fivsev{2.7}

\def\fivzer{2.0}
\def\fivone{2.1}
\def\fivtwo{2.2}
\def\fivthr{2.3}
\def\fivfou{2.4}

\def\sixone{4.1}
\def\sixtwo{4.2}

\def\sevone{5.1}
\def\sevtwo{5.2}
\def\sevthr{5.3}
\def\sevfou{5.4}
\def\sevfiv{5.5}

\def\Sone{\SSS.1}
\def\Stwo{\SSS.2}
\def\Sthi{\SSS.3}
\def\Sfou{\SSS.4}
\def\Sfiv{\SSS.5}
\def\Ssix{\SSS.6}
\def\Ssev{\SSS.7}
\def\Seig{\SSS.8}
\def\Snin{\SSS.9}
\def\Sten{\SSS.10}
\def\Sele{\SSS.11}
\def\Stwe{\SSS.12}
\def\Sthr{\SSS.13}
\def\Sftn{\SSS.14}

\def\Aone{A.1}
\def\Atwo{A.2}
\def\Athr{A.3}
\def\Afou{A.4}
\def\Afiv{A.5}
\def\Asix{A.6}
\def\Asev{A.7}
\def\Aeig{A.8}
\def\Anin{A.9}
\def\Aten{A.8}
\def\Aele{A.11}
\def\Atwe{A.9}

\Title Asymptotic hollowness of lattice simplices

\Abstract An $(n-1)$-tuple $\pmb a:= (a(1), \dots , a(n-1))$ consisting of positive integers is said to be {\it \ah\/} if there exist infinitely many  positive integers $N$ \st  the convex hull in $\R^n = \R^{n\times 1}$,   $K_{\alpha(N)}$, of $\brcs{0,e_1, \dots, e_{n-1},   {\alpha}(N)^T}$ is hollow (no lattice points in its interior), where $e_i$ run over all but the last standard basis element of $\R^n$, and $\alpha(N) = (a(1), a(2), \dots, a(n-1), N)$. Such a tuple is {\it trivial\/} if $\min \brcs{a(i)} =1$. {\par}Nontrivial \ah\ tuples are characterized in terms of modular inequalities,  and turn out to be rare. We show that  for a tuple $\pmb a$, there exists an effectively computable  constant $C$ (depending on $\pmb a$) \st if for  some $N > C$, $K_{\alpha(N)}$  is (not)   hollow, then for all $M > C$, $K_{\alpha(M)}$  is (not)   hollow (respectively); in particular, if $\pmb a$ is asymptotically hollow, the the set of $N$ \st $K_{\alpha(N)}$ is hollow is cofinite in the positive integers. {\par}When $n = 4$, the nontrivial asymptotically hollow triples are completely determined; there are eleven of them, together with   a one-parameter family. 


\noindent {\it David Handelman}%
\plainfootnote{*}{Mathematics Dept, University of Ottawa, Ottawa ON \
K1N 6N5 Canada;  e-mail: dehsg\@uottawa.ca}

{}\plainfootnote{}{AMS(MOS) classification: 11H06, 15A36, 52B20; key words and phrases:  lattice simplex, empty polytope, hollow polytope, asymptotically hollow, facet, width (of a lattice polytope), proscriptive interval}

 \SecT 0 Introduction, definitions, and elementary results 

A {\it lattice polytope\/} is a  subset of Euclidean space  ($\R^n$) that is the convex hull of a finite set of points from $\Z^d$, written $K = \cvx S$ where $S$ is a finite subset of $\Z^d$. We normally also assume that $K$ has interior (that is, contains an open $n$-ball), which is equivalent to the group generated by $S-S$ being of full rank ($n$).
A lattice polytope $K$ is {\it empty\/} if $K \cap \Z^d = \partial_e K$, that is, the only lattice points in $K$ are its extreme points (this is poor terminology, as confusion with an empty set could arise). Finally, a lattice polytope is {\it hollow\/} if it has interior, but no interior lattice points. We are interested in lattice simplices, that is, lattice polytopes that are also simplices. 

Empty and hollow lattice simplices have been heavily studied, especially \wrt classification (up to affine equivalence), e.g.,  [R, W, MMM, IS] and their natural connection to toric varieties, e.g.,  [B, ALPPT], and their connections to other areas, e.g., [B2, H, PR]. This paper is along different lines. A very special class of   families of (mostly) hollow simplices is introduced, and examples provided, with a complete determination in dimension four (in this context, dimension three is not of interest). 

I wanted to find  families of hollow polytopes, requiring only relatively easy computations. This led naturally to the notion of (nontrivial) asymptotic hollowness. It turns out that these are rare (none in dimension $3$, some in dimension four; higher dimensions, mostly unknown)


Section 1  contains some introductory results. Section 2 contains the main result, giving necessary and sufficient conditions to be an asymptotically hollow sequence in terms of modular relations. The short section 3 contains elementary results on facets and width. Section four shows that for every $n > 3$, there exists an asymptotically hollow $(n-1)$-tuple whose least entry is $2^{n-3} (2^{n-2}-1)$; whether this is sharp is unknown. 

Section 5 is a technical section, in preparation for the $n=4$ classification result, although it does apply much more generally. And section 6 contains the classification result for $n = 4$. 
The appendix contains a result that  is used in the determination of \ah\ triples, and might be of interest on its own. 

\SecT 1 Elementary results on lattice simplices

If $A$ is a set, list, or tuple of integers, then we define the {\it content\/} of $A$, $\cont A$, to be the greatest common divisor of the integers in $A$; if $A$ is empty, then $\cont A = 0$. The superscript ${}^T$, applied to matrices, columns, or rows, means {\it transpose.}

Let $d$, $a(i)$, $i = 1,\dots, n-1$ be integers, and let $\aa = (a(1), a(2), \dots, a(n-1),d) \in \Z^{1\times n}$, and assume $d \neq 0$. Form $K_{\aa}$, the convex hull of $\brcs{0, \brcs{e_i}_{i=1}^{n-1}, \aa^T}$, where $e_i$ are the standard basis elements in $\R^n = \R^{n\times 1}$ (columns; the superscripted ${^T}$ indicates transpose). This is a lattice polytope of volume $d$ times that of the standard simplex. If $\cont \aa$ is not $1$, then the edge running from $0$ to $\aa^T$ contains a lattice point. So we normally assume $\cont \aa = 1$, that is, $\gcd \brcs{a(1), \dots, a(n-1); d} = 1$. Now questions arise, which $K_{\aa}$ are empty, and how do we classify them (up to affine equivalence)? 

We form the matrix  $M_{\aa}$ whose first $n-1$ columns are the $e_i$, and the remaining column is $\aa^T$. Elementary row operations on $M_{\aa}$ implement $\GL(d,\Z)$ equivalences between $K_{\aa}$ and the convex hull of the columns of $EM_{\aa}$, where $E$ is the corresponding elementary matrix. By adding or subtracting suitable multiplicities of the last row, we can also assume that $0 \leq a(i) < d$, and in fact, many of the equivalences are really implemented over the finite ring $\Z_d = \Z/d\Z$. Moreover, the presence of a zero as one of the $a(i)$ allows immediate reduction to the smaller size matrix obtained by deleting the $i$th row and column. So we may as well assume $0 < a(i) < d$ (note that $d =1$ yields the standard $n$-simplex, regardless of the $a(i)$). 

For emptiness, we wish to determine when there does not exist a tuple $(\lambda_1, \dots, \lambda_n) $ \st $0 \leq \lambda_i < 1$ for all $i$, $\sum \lambda_i \leq 1$, and $\sum_{i=1}^{n-1} \lambda_i e_i + \lambda_n \aa \in \Z^d$. This amounts to solving an inequality---or more accurately, showing nonexistence of a solution: does there exist an integer  $k$ with $1 \leq k < d$, \st 
$$
 \sum_{\Set{i}{d \text{ does not divide } a(i)k}} \(1 - \Fr{\frac{k a(i)}{d}} \) + \frac kd  \leq 1? \tag 1
$$ 

Here are a couple of elementary classes of examples. 
We normally assume that  {(i)} $1 \leq  a(i) < d$ for all $i$, and  
 (ii) $\cont \aa = 1$

\Lem Lemma \oneone. Sufficient for $K_{\aa}$ to be empty is either of the following conditions. 
\item{(i)} there exists $i$ \st $a(i) = 1$ and $\cont \({\aa}\setminus \brcs {a(i)}\) = 1$; 
\item{(ii)} $\cont (V,d) = 1$ where  
$U = \cup T $ with $T $ varying over  those subsets of $\brcs{1,2,\dots,n}$ \st $d$ divides $\sum_{t\in T} a(t)$,  
 and  $V = \Set{a(i)}{i \in U}$. 

\Rmk These two conditions are very closely related to each other, despite appearances. We will see this when we consider  affine equivalences between various $K_{\aa}$. The second condition is very similar  to the width one property, that will be mentioned later. 

\Pf Suppose there exists a positive integer $k < d$ \st  $\sum \(1 - \Fr{\frac{k a(i)}{d}} \) + \frac kd \leq  1$. 

\noindent (i) If $a(j) = 1$, then $\Fr{a(j)k/d} = k/d$, so $  \(1 - \Fr{\frac{k a(j)}{d}} \) + k/d = 1$. Thus $d$ divides $ka(i)$ for all $i \neq j$. Hence $d$ divides $k\cdot \gcd\brcs{a(i)}_{i\neq j}$; by hypothesis, the second factor is relatively prime to $d$, so $d$ divides $k$, contradicting $k < d$. 

\noindent (ii) If $T$ is a set as indicated, then 
$$\eqalign{
0 &= \Fr{\sum_T a(t)k/d} \equiv \sum_T \Fr{ a(t)k/d} \mod \Z \cr 
& \equiv  \sum_T\(1- \Fr{ a(t)k/d}\) \mod \Z \cr 
}$$
Thus  $\sum_T \(1- \Fr{\sum_T a(t)k/d}\)$ is an integer. Since it is strictly less than $ \sum \(1 - \Fr{\frac{k a(i)}{d}} \) + \frac kd$, which in turn is less than or equal $1$, the sum must be zero; this is impossible, unless it is over an empty set. In particular, this forces $d$ to divide $ka(i)$ for all $i $ in such a $T$, and thus for all $i \in U$. Now the same $\gcd$ argument applies. \qed 

The second criterion is easily refined to the following, by essentially the same proof. 

\Lem Lemma \onetwo. Suppose that (with the notation of the previous lemma), $U$ is nonempty, and $z = \sum \lambda_i e_i + \lambda_n \aa$ is a lattice point  in $K_{\aa} \setminus  \partial_e K_{\aa}$, expressed as convex linear sum of the extreme points other than zero). Then $\lambda_n = k/d$ where $k$ is divisible by $d/\cont \brcs{U,d}$.

Projections onto lower dimensional versions also yield another source of empty polytopes. 

\Lem Lemma \onethr. Suppose that $K_{\aa} \subset \R^n$ is empty, where $\aa = (a(1), a(2), \dots, a(n-1); d)$. If $\aa' = (b(1), b(2), \dots, b(m), a(1), \dots, a(n-1);d)$ and $b = (b(j))$ is any $m$-tuple, then $K_{\aa'} \subset \R^{m+n}$ is empty. \par
If $b$ consists only of zeros, we have essentially a direct sum (in the category of simplices). This suggest a classification of minimal empty simplices, those \st deleting any term causes the resulting lower dimensional simplex to be nonempty.  

\SecT 2 Asymptotic hollowness

Given $\aa(N) := (a(1), a(2), \dots, a(n-1); N) \in \Z^{1 \times n}$ (with $N$
permitted to vary over the positive integers), we wish to determine conditions
guaranteeing or consequences of $K_{\aa(N)}$ being hollow for all sufficiently large $N$, together with upper bounds on the $N$ \st it is not hollow. (Recall that $K_{\alpha(N)}$ is the convex hull---in $\R^{1 \times n}$---of $\brcs{0,e_1, e_2, \dots, e_{n-1}, \aa(N)^T}$.) We say $\pmb a= (a(1), a(2), \dots, a(n-1)) \in \Z^{1 \times (n-1)}$ is {\it \ah} if there exist infinitely many positive integers $N$ \st $K_{\alpha(N)}$ is hollow. A simple argument shows that  if   there exists $i$ \st $a(i)=1$, then $K_{\alpha(N)}$ is hollow; so we say an \ah\ tuple $\pmb a$  is {\it nontrivial\/} if $a(i) \neq 1$ for all $i$. 

It turns out  that this notion is not interesting if $n = 3$; that is, a pair of positive integers is \ah\ only if one of its entries is $1$.   This follows from  the description empty  three-dimensional lattice simplices ([R, W]), but we can give a direct and elementary argument. 

\Lem Lemma \fivzer. Let $(a,x)$ be an ordered pair of positive integers, with $2 \leq a \leq x$.  Let
 $N$ be an integer such that $(N-x)(a-1) \geq x^2$ (that is, $N > x + x^2/(a-1)$). On defining $\alpha = (a,x,N)$, we have that  $\cvx \{ 0, e_1, e_2, \alpha^T\}$  has an
interior lattice point. 

\Pf We will find $\lambda_i > 0$ such that $\lambda_1 e_1 + \lambda_2 e_2 +
\lambda_3 \alpha^T$ is a lattice point and $\lambda _1 + \lambda_2 + \lambda_3
\leq 1$. Write $N = mx + r$ where $1 \leq r \leq x$ (note: not $0 \leq r \leq x-1$), and set $\lambda_3 = m/N$.
As $a \leq x$ and $mx < N$, we have $ma < N$. Then we set $\lambda _1 = 1-ma/N >
0$, and $\lambda_2 = 1-mx/N = r/N > 0$. It is clear that $\lambda_1 e_1 +
\lambda_2 e_2 + \lambda_3 \alpha^T$ is the  lattice point $(1,1,m)^T$. Now $\lambda _1 +
\lambda_2 + \lambda_3  = 1-ma/N + r/N + m/N \leq 1+ (-ma + x + m)/N$, so this
will be at most one if $m(a-1) \geq x$. But $m = (N-r)/x \geq (N-x)/x$, and our
assumption was that $((N-x)/x) (a-1) \geq x$. \qed

In the case that $n =4$ (that is, triples), nontrivial \ah\ triples exist, and are completely described in section \SSS. 

We introduce the following notation (this is open to suggestions for improvement). Let $x$ and $y$ be integers with $x > 1$. Define 
$$
\brm x,y. = \cases
r  & \text{if $1 \leq r \leq x-1$, and $r$ is the remainder of $y$ modulo $x$}\\
x & \text{if $x$ divides $y$}.\\
\endcases$$
This is the usual remainder function, except that instead of having zero as a value, it has $x$ as a value. Obviously, $\brm x,y. \equiv y \pmod x$. To avoid excessive parenthesization, we often write $\brm x,(ty).$ simply as $\brm x,ty.$.

\Lem Theorem \fivone.  (a) If $N \geq \max_{i\neq j}\brcs{(a(i)-1)a(j)}$, then sufficient for
  $K_{\aa(N)}$ to be hollow is that  for all $i,t$, 
$$
\sum_{j \neq i}  \brm a(i),ta(j). \leq t + (n-3)a(i).
$$ 
\item{(b)} If additionally $N > \( \sum a(i) - 1\)\cdot \max_i \brcs{a(i)-1}$,
then this condition is also necessary.  

\Pf Let $b$ be a lattice point in $K_{\aa(N)}$ other than an extreme point. Then
there exists a positive integer $k <N$ \st 
$$\eqalign{ 
b &=  \sum_{\Set{i}{N \text{ does not divide } k a(i)}} \(1 - \brcs{\frac{k
a(i)}{N}} \)e_i+ \frac kN   \aa(N)^T \cr 
 1 & \geq  \sum_{\Set{i}{N \text{ does not divide } k a(i)}} \(1 - \brcs{\frac{k
a(i)}{N}} \)+ \frac kN . 
}\tag 1$$
We can rewrite the inequality as 
$$ 
 \sum_{\Set{i}{N \text{ does not divide } k a(i)}}  \brcs{\frac{k a(i)}{N}} -
\frac kN \geq \left| \Set{i}{N \text{ does not divide } k a(i)}\right|  - 1.
\tag {$2'$} 
$$ 
Since $K_{\aa(N)^T}$ is a simplex, the decomposition will be unique (given by
the barycentric coordinates from the extreme points, $\brcs{0; e_1,\dots, e_{n-1}; \aa(N)^T}$).
Then $b$ will be in the interior if and only if there does not exist $i$ \st $N$
divides $ka(i)$ and the inequality in the second line is strict. So let us
assume that $b$ is in the interior. In that case, $ \Set{i}{N \text{ does
not divide } k a(i)}  $ is just the entire set of $i$s. So assuming $b$ is an interior point, the inequality becomes, 
$$ 
 \sum_{i=1}^{n-1}  \brcs{\frac{k a(i)}{N}} -
\frac kN \geq n-2.
\tag {$2$} 
$$
And of course, (2) holding for some $k $ is both necessary and sufficient for $K_{\aa(N)}$ to have an interior lattice  point.

Now consider the finite set of nonzero rational numbers, 
$$
U= \Set{\frac{t}{a(i)}}{i = 1, 2, \dots, n-1; \text{ and } 1 \leq t \leq a(i)-1}
$$ 

Then $U^c = (0,1) \setminus U  $ is a finite union of open intervals, whose
endpoints are of the form $t/a(i)$. It can happen that $t_1/a(i_1) = t_2/
a(i_2)$ with $i_1 \neq i_2$, but this is not a problem (and will not occur if
the  $a(i)$ are pairwise relatively prime). 

Now let $b$ be an interior point, with corresponding $k$. Since $N$ does not
divide $ka(i)$ for any $i$, it follows that $k/N $ does not belong to $U$. Hence
there exists a unique open interval $I= (t_1/a(i_1), (t_2/a(i_2))$ that is a
component of $U^c$ and contains $k/N$. From the definitions of $U$, it follows
that $I$ contains no elements of the form $t/a(i)$. Thus for each $i$, the
function $x \mapsto \flo {x a(i)}$ is constant on $I$ ($x$ is a real variable,
the continuous version of  $k/N$). Thus the function $I \to \R$ given by 
$$
x \mapsto \sum \brcs{ xa(i)} - x 
$$
 can be rewritten in the form
$$
F:x \mapsto \sum  xa(i) - x  - C
$$ 
where $C$  is the constant $\sum \flo{xa(i)}$.  The derivative is thus $\sum
a(i) -1$---which is a positive integer. In particular, $F$ is strictly increasing on $I$.
Hence if $k' > k$ is another integer with $k'/N \in I$, then $F(k'/N) > F(k/N)$.
But from (2), this implies that $F(k'/N) > F(k/N) > n-2$, and thus $k'/N$ yields
another interior point. 

The largest $k'$ \st $k'/N \in I$ is given by $k' = \flo {N t_2/a(i_2)}$
(obviously $k' < Nt_2/a(i_2)$). Hence given a solution (an interior lattice
point) parameterized by $k$, we obtain a solution arising from $k' = \flo
{Nt_2/a(i_2)}$. In particular, to decide whether $K_{\aa(N)}$ is hollow, we only
have to investigate the choices for $k$ of the form $\flo{Nt/a(i)}$ as $i$
varies over $1,2,\dots, n-1$ and $1 \leq t < a(i)-1$. The number of such choices
is $\sum a(i)-1$ (under some conditions, we can halve this; see Lemma \fivthr).
The number of trials is independent of $N$, except when the latter is relatively
small. We now investigate what conditions this imposes. Relabel $k' $ as $k$. 

Let  $0 < \rho_{i_2,t} \equiv \rho < a(i_2)$ be the remainder of $Nt_2$ modulo
$a(i_2)$, i.e., $k/N = (t_2 - \rho/N)/a(i_2)$ (we are assuming $N > \max a(i)$, as
usual). 

\comment
Also, for $j \neq i_2$, us the notation $\alpha_{i_2,t}(j) \equiv \alpha(j)$ to
be the remainder of  $t a(j)$ modulo $a(i_2)$, except that we only permit values
of $\alpha(j)$ in the interval $\brcs{1,2,\dots , a(i_2)}$ (so if $a(i_2)$
divides $t a(j)$, the corresponding $\alpha(j) = a(i_2)$). 
\endcomment

We have $k a(i_2) = Nt_2 - \rho$, so $\brcs{ka(i_2)/N} = 1 - \rho/N$. For $j
\neq i_2$, we have  
$$\eqalign{
\frac{ka(j)}{N}  & = \frac{(Nt- \rho)a(j)}{a(i_2)N} \cr & = \frac{ta(j)}{a(i_2)}  -
\frac{\rho a(j)}{N a(i_2)} \cr 
& = \frac{\brm a(i_2),ta(j).}{a(i_2)} + \flo{\frac{ta(j)-1}{a(i_2)}} - \frac{\rho a(j)}{N
a(i_2)}, \quad \text{and thus}\cr 
\brcs{\frac{ka(j)}{N} } & = \brcs{\frac{\brm a(i_2),ta(j). -  a(j)\rho/N}{a(i_2)}}.\cr 
}$$

Now we make assumption (m): 
$$
N \geq \frac{\rho_{i,t}a(j)}{\brm a(i),ta(j).}\quad \text{for all $i$ and $t < a(i)$}.
\tag m   $$ (This occurs if  $N > a(i)a(j)$ for all $i \neq j$; we will slightly refine this later.) Under this assumption,  
$$
\brcs{\frac{ka(j)}{N} }  = \frac{\brm a(i),ta(j). -  a(j)\rho/N}{a(i)} . 
$$ 

Plugging everything into the left side of (2), 
$$\eqalign{ 
\sum \brcs{\frac {ka(i)}{N}}  -\frac kN & = 1 +  \frac 1{a(i_2)} \(\sum_{j \neq
i_2} \brm a(i),ta(j). -t - \frac{\rho}N\(\sum_{i=1}^{n-1}  a(j) - 1\)\);   \text{ this exceeds $n-2$ iff} \cr 
\sum_{j \neq i_2} \brm a(i),ta(j). -t & >  \frac{\rho}N\(\sum_{i=1}^{n-1}  a(j) - 1\) +
(n-3) a(i_2).
}$$

Relabelling the terms, this means that in the presence of (m),  necessary and sufficient for
$K_{\aa(N)}$ to be hollow is  that for all $i,t$,
$$ 
\sum_{j\neq i} \brm a(i),ta(j).  - t - \frac {\rho_{i,t}}N  \( \sum_1^n a(j) - 1\) \leq
(n-3)a(i).
$$ 
 In particular, 
$$ 
\sum \brm a(i),ta(j).  - t  \leq  (n-3)a(i) \tag 3 
$$ is sufficient (in the presence of (m)) for $K_{\aa(N)}$ to be hollow. 

If we make assumption (M), 
$$
 \frac {\rho_{i,t}}N  \( \sum_1^n a(j) - 1\) < 1,  \tag M
$$
 then  (3) is also necessary. 
Sufficient for (M) to hold is 
$$
N >\( \sum a(j) - 1\) \rho_{i,t},  
$$ 
or a little less generally (but subscript-independent), $N > \( \sum a(i) -
1\)\cdot \max_i \brcs{a(i)-1}$. 
Similarly, $N \geq \max_{i\neq j}\brcs{(a(i)-1)a(j)}$ is sufficient for (m) to hold. \qed

For example, it is easy to verify from these criteria that $\alpha(N): = (3,5,7;
N)$ is hollow for $N \geq 30$. Verification by hand reveals all other choices of
$N$ also yield hollow polytopes. Examination of the faces reveals that provided
$\gcd (N,7) = 1$, all the facets are themselves empty, and thus we obtain that
if $7$ does not divide $N$, then $K_{\aa(N)}$  is empty (when $7$ divides $N$,
two of the facets are not empty). Here $n = 4$ is small enough that  emptiness
of the facets is easy to check. In general, checking that the facets are
themselves empty (or at least relatively hollow), involves reduction of the $N$
to much smaller numbers---specifically $\gcd\brcs{a(i),N} $ for each $i$, and
$\gcd \brcs{\sum a(j)-1,N} $---so computations are going to be much less 
cumbersome. 

In particular, a sufficient condition for emptiness is that $K_{\aa(N)}$ be
hollow, and for all $i$,  $\gcd\brcs{a(i),N} = 1 = \gcd \brcs{\sum a(j)-1,N}$.

\Lem Corollary \fivtwo. Let  $\aa(N) = (a(1), \dots, a(n-1);N)$ where $0 < a(i) $, and suppose there exists an integer $M > \( \sum a(i) - 1\)\cdot \max_i \brcs{a(i)-1} $ \st $K_{\aa(M)}$ is hollow (not hollow, respectively). Then $K_{\aa(M)}$ is hollow (not hollow, respectively) for all $N > \( \sum a(i) - 1\)\cdot \max_i \brcs{a(i)-1} $.

We can improve the result of Theorem \fivone, in that we can restrict to $1 \leq t \leq a(i)/2$ rather than all $1\leq t \leq a(i)$.

\Lem Lemma \fivthr. Assume that condition (m) applies to $\alpha \equiv \alpha(N)$. Suppose that for $k = \flo{Nt/a(i)}$ with $a(i)/2 < t < a(i)$, the corresponding point, 
$$
v:= \sum_{\Set{j}{a(i) \text{ does not divide } N a(j)}} \(1- \brcs{\frac{ka(j)}{N}} \)e_i + \frac kN  \aa^T \in K_{\aa}.
$$  
Then the same is true with $k$ replaced by $k_0 = \flo{Nt_0/a(i)}$, where $t_0 = 2 t - a(i)$. If the former is in the interior, then so is lattice point arising from $k_0$.

\Pf Since $K_{\aa}$  is a simplex, the displayed condition is simply,
$$
\sum_{\Set{j}{a(i) \text{ does not divide }  
N a(j)}} \(1- \brcs{\frac{ka(j)}{N}} \) + \frac kN  \leq 1.  
$$
In the presence of (m), this boils down to 
$$
\sum_{j \neq i} \brm a(i),ta(j). - t -\frac {\rho}{N}\(\sum_{j=1}^n a(j) -1\) \geq (n-3)a(i) ,
$$
where $\rho = \rho_{i,t}$
. Define $\alpha_2 = \brm a(i),t_0 a(j). $ and $\rho_2 = \rho_{i,t_0}$. Since the right hand side differs from $(n-2)a(i)$ by less than $a(i)/2$, and there are only $n-2$ positive terms on the left each of which has maximum value at most $a(i)$, it follows that each $\alpha(j) > a(i)/2$, and thus $\alpha_2(j) = 2\alpha (i) - a(i)$. Also $\rho_2 $ is either $2\rho$ (when $\rho < a(i)/2$) or $2\rho - a(i)$; write $\varepsilon $ to be $1$ if $\rho/a(i) > 1/2$, else zero.  Then $\rho_2 = 2\rho -\varepsilon$. The corresponding expression on the left for $k_0$ is 
$$\eqalign{ 
2 \sum \brm a(i),ta(j).  -& (n-2)a(i) - (2t-a(i)) - \frac{2\rho - \varepsilon}{N}\(\sum_{j=1}^n a(j) -1\)\cr 
& = 2\(\sum \brm a(i),ta(j).  - t -  \frac{\rho}{N}\(\sum_{j=1}^n a(j) -1\) \)-(n-3)a(i) + \frac{\varepsilon a(i)}{N}\cr 
& \geq  (2n-6 - (n-3))a(i)\cr  & = (n-3)a(i). 
}$$
So the lattice point associated to $k_0$ belongs to $K_{\aa}$. If $v$ belongs to the interior, then all the coefficients in the barycentric decomposition are nonzero, and their sum is not one. It is easy to check that this will also apply to the new lattice point. \qed

\Lem Proposition \fivfou. Suppose that $K_{\aa}$ admits a non-extreme (interior) lattice point, and that condition (m) holds. Then there exist positive integers  $i$ and  $s$ with $ s < a(i)/2$ \st $k = \flo{Ns/a(i)}$ yields a lattice point (in the interior of) $K_{\aa}$.  

\Rmk Sufficient for condition (m) to hold is that $N > \max_{i \neq j} a(j)(a(i) - 1)$. 

\Pf The function $f:x \mapsto 2x - 1$ has a fixed point at $x =1$ and for $0 < x < 1$, the sequence of iterates $(f^n(x))$  converges to zero. Hence given $x = t/a(i)> 1/2$, there exists a smallest integer $u$ \st $f^{u}(x) < 1/2$. Now the lemma above applies.  \qed

Can we do better, that is, can we drastically weaken the condition (m) hypothesis to obtain a result on minimal $t$? The proof of Lemma \fivfou\  uses what appear to be flabby estimates. To this end, define $S_+ = \Set{j}{\alpha(j) > \rho a(j)/N}$, and $S_{-}$ its complement in $\brcs{i}^c$. If any of $|S_{-}| \leq (n-2)/2$, or $ \brcs{Nt/a(i)} \leq 1/2$, or $|S_{-}| = |S_{2,-}|$  (the latter is the corresponding set for $t_0$), then the first round of iteration succeeds, that is, the lemma applies---once. Unfortunately, none of these conditions are likely to be preserved in the next iteration, so we cannot continue the argument. (A notable exception occurs if $t/a(i) < 3/4$, because only one iteration is required.) A preferable condition would be something like $N > 2\max a(j)$.

 The upshot is that at least if (m) applies, then we only have test to the cases with $t < a(i)/2$. So at most $\sum \flo{a(i)/2}$ congruence inequalities have to be checked to decide hollowness.

\Lem Corollary \fivfiv. If $N \geq \max_{i\neq j}\brcs{(a(i)-1)a(j)}$, then sufficient for
  $K_{\aa(N)}$ to be hollow is that  for all $i,t$ with $t < a(i)/2$, 
$$
\sum_{j \neq i} \brm a(i),ta(j).\leq t + (n-3)a(i).
$$
In particular,  the displayed condition is necessary and sufficient for $\pmb a = (a(1), \dots, a(n-1))$ to be asymptotically hollow. 
\qed

The following is elementary, but very useful  in reducing the number of times we have to decide whether the inequality in Corollary \fivfiv\ holds. 

 \Lem Lemma \fivsix. Let $\pmb a = (a(1), \dots, a(n-1)) \in \Z^{n-1}$ with all $a(i)\geq 2$. Suppose there exists $j$ and $S \subset \brcs{1,2,\dots, n-1}\setminus \brcs{j}$ \st $\sum_{s\in S} a(s) \equiv r \pmod {a(j)}$.
Either of the following conditions is sufficient to guarantee that  for all $t\in \N$, $\sum_{i \neq j}  \brm a(j),a(i)t. \leq t + (n-3)a(j)$.{ }
\item{(a)}  $r = 1$.  { }
\item{(b)}  $r = 0$ and there exists $s \in S$ \st $a(s)$ is not divisible by $a(j)$.

\Pf   Define $s_t = \sum_S \brm a(j),a(i)t .\equiv t \pmod {a(j)}$; then $s_t \equiv tr \pmod {a(j)}$. It is sufficient to verify the condition for $t < a(j)/2$, by Corollary \fivfiv. 

\noindent  (a) We have $s_t=\sum_S \brm a(j),a(i)t. \equiv t \pmod a(j)$, so we can write $s_t - t = ka(j)$ for some integer $k$. Since $s_t \leq |S|a(j)$, it follows that $k \leq |S|-1$. Then  
$$
\eqalign{\sum_{i \neq j}  \brm a(j),a(i)t. - t &\leq \(s_t - t \) + \sum_{S^c \setminus\brcs{j}} \brm a(j),a(i)t. \cr  
& \leq \(|S|-1\)a(j) + |S^c \setminus \brcs{j}| a(j) \cr 
& = ((n-1) - 2 ) a(j).\cr
} $$

\noindent (b) We can write $s_t = ka(j)$, where $k \leq |S|$. However, not all $a(s)$ are divisible by $a(j)$, hence for fixed not all $\brm a(j),a(s)t. = a(j)$; hence $k < |S|$. Now the same argument applies, and we deduce (in this case) that $\sum_{i \neq j}  \brm a(j),a(i)t. \leq (n-3)a(j)$ (with no $t$ term on the right). \qed 

This generalizes; we have a result for specific $t$. The proof is essentially the same. 

\Lem Lemma \fivsev. Let $(a(i)) \in \Z^{n-1}$. Suppose that for some $j$ and positive integer $t$, there exists a subset $S \subset \brcs{1,2,\dots,n-1}\setminus \brcs{j}$, \st 
$$
t\sum_{i \in S} a(i)  \equiv z \mod {a(j)},
$$
for some $z$ with $1 \leq z \leq t$, or if $z = 0$  and there exists $i \in S$ \st $a(i)$ does not divide $a(j)$. Then 
$$
\sum_{i \neq j} \brm a(j),a(i)t. \leq t + (n-3)a(j).
$$

\Pf Since $\sum_{i \neq j} \brm a(j),a(i)t. \equiv t\sum_{i \in S} a(i)$, we have that $\sum_{i \neq j} \brm a(j),a(i)t. \leq (|S|-1) a(j) + z$, and this is bounded above by $(|S|-1) a(j) + t$. Now the same finishing  arguments as those in Lemma \fivsix(a,b) complete the proof. 
\qed 

For example, suppose we have the quadruple $(3,29, 38, 66)$ ($n=5$), and wish to determine whether the  inequalities  modulo (say) $29$ hold. First, reduce the three remaining numbers modulo $29$, so we are dealing with $3,9, 8$. Using $S = \brcs{3}$, applying Lemma \fivsev, we have that the inequality holds for $10 \leq t \leq 14$. From $S = \brcs{3,8,9}$, we see that it holds for $t = 3, 6, 9$; with $S = \brcs{8}$, we obtain $t = 4,8$. When $S = \brcs{3,9}$, the inequality holds for $t = 5$, and when $S = \brcs{8,9}$, $t = 7$ works.  Since the inequality trivially holds when $t = 1$ or $2$, it holds for all $t$. In fact, this quadruple is \ah, as we will see in Example \sevfou.

\SecT 3 Facets and width 

Let $A \subset \Z^n$, and form the group generated by $A-A$ (the set of differences of elements of $A$), $H(A)$. If the torsion subgroup of  $\Z^n/H(A)$ has order $m$ (that is, the torsion subgroup of $\Z^n/H(A)$ has $m$ elements---it will always be finite), we say $A$ generates a subgroup of cotorsion $m$.

Now let $K$ be a lattice polytope with interior, and $F$ a facet of $K$. Set $A = F \cap \Z^d$. Then it is easy to see that the torsion-free rank of $H(A)$ is $n-1$ (since we assumed at the outset that $K$ has interior). Recall the notation: if $\alpha = (a(1), \dots, a(n-1),d) \in \Z^{1 \times n}$, then $K_{\alpha}$ is the convex hull of $\brcs{0, e_1, \dots , e_{ n-1}, \alpha^T}$, and this has interior if and only if   $d$ is not zero---and we normally assume that $d \geq 1$. 

Define the {\it volume\/} of $F$ to be the $m$ \st  $A$ generates a subgroup (of $\Z^n$) of cotorsion $m$. If $m =1$, then $F$ is a {\it standard\/} (or {\it unimodular\/}§) ($d-1$) simplex, and in general $m$ is the normalized volume of $F$. 

Now we look at the facets of $K_{\aa}$. There are of course $n+1$ of them. One of them is the standard $d-1$ simplex, obtained by  deleting $\aa$ from the list of vertices. Deleting one of the $e_i$ yields a facet of (normalized) volume $\gcd\brcs{a(i),d}$, as is easy to verify. Deleting $0$ is slightly more complicated, in that we have to translate $\cvx \brcs{\brcs{e_i}_{i=1}^n \cup \brcs{\aa^T}}$ so as to hit the origin, e.g., by subtracting $e_1$ from everything, and it is routine that the volume of this face is $\gcd\brcs{\sum a(i) -1, d}$. We thus have a list of the volumes of the facets of $K_{\aa}$: 
$$ 
1,\gcd (a(1),d), \gcd(a(2),d), \dots, \gcd(a(n-1),d);\gcd \(\sum a(i) -1, d\), 
$$
the last term corresponding to the facet obtained by deleting $0$.

\Lem Lemma \throne. Let $K = K_{ \aa}$. \item{(i)} The number of facets of  $K$ that  are standard simplices is $1$ plus the number of $j = 1, 2, \dots ,n$ in the augmented row $b = (b(j))$ (with $b (j) = a(j)$ if $j \leq n-1$, and $b(n) = \sum a(i) -1$) \st $\gcd \brcs{b(i),d} = 1$; this is an affine integral invariant of $K$. \item{(ii)} All of the facets are standard simplices if and only if  $\gcd\( \sum a(i) -1, d\) = 1$ and $\gcd (a(i),1) = 1$ for all $i \leq n-1$. 

When $d$ is prime and all of the $a(i)$ are not divisible by $d$, then there are either $n-1$ or $n$ facets that are standard simplices, depending on whether $\sum a(i) -1$ is divisible by $d$. Both situations can occur, and constructions are quite easy, even if we also require emptiness. 

If $K_{\alpha}$ is hollow, then  necessary and sufficient that it be empty is that every  facet is hollow (relative to the corresponding sublattice); in particular, this occurs if all the  facets are standard.

\noindent {\it Width } Let $K$ be a lattice polytope, and let $\phi \in \Z^{1\times n}$ (equivalently, $\phi$ is a homomorphism $\Z^d \to \Z$). Define the {\it width relative to $\phi$ of $K$}, $\wid_K (\phi)$ to be $\max \phi|K - \min \phi|K$. Since linear functionals take their extreme values on extreme points (which here are lattice points), this is an integer, and if we take the infimum of $\wid_K (\phi)$ over all nonzero $\phi$, we see that the infimum is achieved. This defines $\wid (K)$. Obviously, this is an $\agl\Z$ invariant, and is never zero.

For $K = K_{\aa}$, width one is easy to characterize. 

\Lem Lemma \thrtwo. $K_{\aa}$ has width one iff either there exists a subset $T \subset \brcs{1,2,\dots, n-1}$ \st $\sum_{i \in T} a(i) \equiv 0 \pmod d$ or $\sum_T a(i) \equiv 1 \pmod d$.

\Pf We can obviously assume $d > 1$. First suppose that $K_{\aa}$ has width one. Since $0 \in K_{\aa}$, the integer $0$ is a value of any $\phi \in \Z^{1 \times n}$. If $\phi$ achieves the minimal width, then either all of its other values are either $1$, or all of its other values are either $-1$. Replacing $\phi$ by $-\phi$ if necessary, we can assume that its values on $K$ are all in $\brcs{0,1}$, and since $K$ contains interior and $\phi$ is supposed to be nonzero, $1$ is a value. 

Writing $\phi $ as a row, from the values on the first $n-1$ of the $e_i$, we see that $\phi = (\epsilon(i); x)$ where $\epsilon (i) \in \brcs{0,1}$ for all $i \leq n-1$.  Set $T = \Set{i}{\epsilon(i) = 1}$. Then $\phi \cdot \aa = \sum_T a(t) + xd$. If this is  zero, then we are done; otherwise, it equals $1$, and again we are done. 

Conversely, assume such a subset $T$ exists.  Write $\sum_T a(T) = \epsilon + md$ where $\epsilon \in \brcs{0,1}$, and define $\phi = (\epsilon(1), \epsilon(2), \dots, \epsilon(n-1), -m)$ where $\epsilon(i) = 1$ if $i \in T$ and $0$ otherwise. \qed 

This extends, somewhat. 

\Lem Lemma \thrthr. Let $K = K_{\aa}$, and let $U$ be the set of integers that appear in any of the augmented forms resulting $(a(i))$ and its iterates under the actions described in section\,2. Rewrite the elements  of $U$ modulo $d$ so that $-d/2 < u \leq d/2$ for all $u \in U$, and suppose that $0$ never appears (that is, no edges contain lattice points other than the extreme points). Set $s = \min \brcs{\Set{u}{u \in U; u> 0} \cup  \Set{1+ |u|}{u \in U; u< 0}}$. Then $\wid (K) \leq s$. 

\Pf Pick the corresponding $u$ to achieve $s$. There exists an equivalent simplex $K'= K_{\aa'}$ \st $u +md$ is the $i$th entry  of $\aa'$ for some $i \leq d-1$. Then set $\phi $ to be the row with a $1$  in the $i$ the entry, $-m$ in the $d$th entry, and zeros elsewhere. Then the values of $\phi$ on the extreme points are $0, 1$, and $u$. If $u < 0$, then $\wid_{K'}(\phi)$ is $1 - u = 1 + |u| = s$, whereas if $u > 0$, then $\wid_{K'}(\phi) = u = s$. Thus $\wid (K') \leq s$. The result follows from the observation that the width is an affine invariant. \qed 

Unfortunately, the inequality can be strict. Take $\aa = (2,3,3,3,4; 12)$. Then because there are no equivalences arising from inverting modulo $d$, there is only one augmented row to check,  $(-2,2,3,3,3,4)$. It contains no ones, so the $s$ in this case is $2$. However, $2 + 3 + 3 + 4 = 12$, so by Lemma \thrtwo, the width is one. This can be fixed by restricting to projections onto lower dimensional spaces, e.g., $(2,3,3,4; 12)$ obtained by deleting one of the terms is still empty, only this time $1$ appears in the augmented row $(1,2,3,3,4)$ (and there is only one other augmented row, $(1,-2,-3,-3,-4$)). So we could take the infimum of the $s$ obtained by deleting some of the terms.

\input addonstoemptiness

\long\def\Rf[#1] #2, #3. #4\par%
{\vskip 2pt \itemitem{[#1]} #2, {\it #3,} #4\par\vskip2pt}

\SecT References

\Rf [ALPPT] {A Atanasov, C Lopez, A Perry, N Proudfoot, M Thaddeus}, Resolving toric varieties with Nash blow-ups. Experimental Mathematics 20 (2011) 288--303.

\Rf [B] A Borisov, On classification of toric singularities. J Math Sci NY 94 (1999) 1111--1113. 

\Rf [B2] ---\!---\!---, Quotient singularities, integer ratios of factorials, and the Riemann Hypothesis. 

\Rf [H] D Handelman, Invariants for critical dimension groups and permutation-Hermite equivalence, M\"unster J of Math. 11 (2018) 49--156.

\Rf [IS] O Iglesias-Vali–o and F Santos, Classification of empty lattice 4-simplices of width larger than two. 

\Rf [MMM] S Mori{, DR Morrison, and I Morrison}, On four-dimensional terminal quotient singularities. Math Comput 51 (1988) 769--786.

\Rf [PR] V Powers and B Reznick, A note on mediated simplices. J Pure and App algebra, 225 (2021)

\Rf [R] JE Reeve, On the volume of lattice polyhedra. Proc London Math Soc 7 (1957) 378--395 

\Rf [W] GH White, Lattice tetrahedra. Can J Math (1964) 389--396.

{}

\vskip 10pt

Mathematics Department, University of Ottawa, Ottawa ON  K1N 6N5, Canada; dehsg\@uottawa.ca

\end

general stuff (if not already done)

\Lem Lemma. Suppose for $\pmb a = (a(1), \dots, a(n-2)) \in \N^{n-2}$, at least one of the proscriptive intervals  $[a(i)/j, s/f(i,j))$ is nonempty [here $s = \sum a(i) -1$. Then there exist at most finitely many positive integers $y$ \st $\pmb a_y:= (a(1), \dots, a(n-2), y)$ is \ah. 

\Pf From the general construction of proscriptive intervals, $y \notin \cup_{t\in \N} [ta(i)/j, ts/f(i,j))$. As $[a(i)/j, s/f(i,j))$ is nonempty by hypothesis, the union contains a ray, that is, $[N,\infty)$ for some positive integer $N$. Hence $y \leq N-1$ is bounded. \qed For this result, we do not require the $n-2$- or $n-1$-tuples to be in increasing order. So a statement equivalent to the following is that if there is a proper subsequence \st for infinitely many positive integer, adjoining that integer yields an \ah\ sequence, then the subsequence is \ah. 

\Lem Proposition. Let $\pmb a = (a(i))_{i=1}^{n-2}$. If there exist infinitely many $y$ \st $\pmb a_y:= (a(1), \dots, a(n-1),y)$ is \ah, then $\pmb a$ is already \ah.

\Pf  By the previous result, all the proscriptive intervals $ [a(i)/j, s/f(i,j))$ are empty. It easily follows that for any positive integer $y$, $\sum_i \brm y,a(i)t. \leq (n-3)y + t$. In particular, if we set $y = a(i)$, we deduce (since $\brm a(i),a(i)t. = a(i)$) that $\sum_{j\neq i}\brm a(i),a(j)t. \leq (n-4)a(i) + t$. But this (for all $i = 1,2,\dots, n-2$) is precisely the set of conditions equivalent to $\pmb a$ being  \ah. 

Quantitative proof of preceding 

\Lem Proposition. Suppose that $\pmb a':= (a(1), a(2), \dots, a(n-2)) \in \N^{n-2}$ is not \ah.  \item{(a)} There exist positive integers $i$ and $m$ with $i \leq n-2$ and $m\leq a(i)/2$ \st $$
\sum_{l\neq i} \brm a(i), a(l)m. > (n-4) a(i) +m ,\tag1 $$ and the corresponding proscriptive interval $[a(i)/m, s/H(i,m))$ (where $s = \sum_{j=1}^{n-2} a(i) - 1$) is nonempty. 
 \item{(b)} Suppose that $i,m$ are positive integers \st (1) holds.  If $\pmb a:= (a(1), \dots, a(n-2),y) \in \N^{n-1}$ is \ah\ for some positive integer $y$, then  
$$ y  < \frac{a(i)\(1 + s-m}{2m}\). $$

\Pf The first part of (a) is a consequence of the general criterion for asymptotic hollowness, applied with $n-1$ replacing $n-2$. Consequent nonemptiness of the interval was shown earlier.  \noindent (b) Set $I = [a(i)/m, s/H(i,m))$, so that $y \not\in \cup_{t\in \N}tI$ (by earlier computation). There exists $t_0 \in \N$ \st for all $t \geq t_0$, we have $st/H(i,m) \geq a(i)/m$, estimated as follows. The inequality can be rewritten $$\eqalign{
t_0\(\frac{ms}{H(i,m)} - a(i) \)& \geq a(i); \quad{\text{so}}\cr 
t_0 & \geq \frac{a(i)H(i,m)}{ms - a(i)H(i,m)} \cr 
&= \frac{a(i)(n+m-4) + m\sum_{l\neq i} a(l) - \sum_{l\neq i}\brm a(i),ma(l). }{ms - \(a(i)(n+m-4) + m\sum_{l\neq i} a(l) - \sum_{l\neq i}\brm a(i),ma(l).\)}\cr
 &=  \frac{a(i)(n-4) + m + ms  - \sum_{l\neq i}\brm a(i),ma(l). }{m- a(i)(n-4)  +\sum_{l\neq i}\brm a(i),ma(l).}\cr &= \frac{\frac{a(i)(n-4)}m + s+1  - \frac{\sum_{l\neq i}\brm a(i),ma(l).}m }{1- \frac{a(i)(n-4)}m   +\frac{\sum_{l\neq i}\brm a(i),ma(l).\)}m}\cr 
}$$

If $t_0$ is set to the ceiling of the last expression, than any larger $t$ will also satisfy the inequality, and thus the ray $[a(i)t_0/m, \infty)$ is contained in $\cup_{t\in \N} tI$. In particular, $y \leq a(i)t_0/m$. The denominator is at least as large as $2$ (by (1)), and the numerator is bounded above by $s+1 - m$. So we can take $t_0 = \ceil{(s+1 - m)/2}$, e.g., $t_0 \leq (s-m)/2 + 1$. Thus $a(i)t_0/m \leq (1 + (s-m))a(i)/2m$.  \qed \par \par If there are numerous pairs $(i,m)$ \st (1), then the the upper bound can be made much smaller (often), since there will be overlapping of multiples of different intervals. 

\noindent $\aa= (3, 29, 38, 66)$. Modulo either $d = 29$ or $38$, the lemma applies, so we need only test for $t \leq d/3$.  more needed 
 
   \noindent $\aa= (11, 29, 38, 66)$. Since $29 + 38 \equiv 1 \pmod {11}$ and $11 + 66 \equiv 1 \pmod {38}$, we need only consider the inequalities  modulo $d = 29$ and $t \leq 14$. The remainders are $11, 9, 8$. In this case, we can use method of proof of the lemma. Set $x = \brm 29,8t.$ so that $\brm 29,9t. \equiv x + t$ and $\brm 29,11t. \equiv x + 3t$.  Moreover, since $29$ is prime, $x \leq 28$. 

 If $x + 3t \leq 29$, then the inequality yields 
$$
 x + 3t + x + t + x > t + 58.
 $$ 
 Thus $3x + 3t > 58$, so $x+t > 19\frac13$, that is $x + t \geq 20$. Since $x + 3t \leq 20$, we deduce $2t \leq 9$, so $t \leq 4$. For the resulting four calculations, it is easy to check that that the inequalities hold. 
 
 Now suppose that $29 - 3t < x \leq 29- t$. Then $\brm d,11t. = x + 3t - 29$, and the resulting inequality is $x + 3t - 29 + x+t + x > t + 58$, or $3(x+t) > 87$, contradicting $x + t \leq 29$. 

 Finally suppose that $ x > 29- t$. Then either $29 < x + 3t \leq 58$ or $58 < x+3t \leq 87$, so that either $\brm 29,11t. = x+3t - 29$ or $x + 3t-58$. Either way, we obtain a something; in the first case, $x+ 3t - 29 + x+t - 29 + x > t + 58$ yields $x > 38\frac{2}3 - t$, forcing $t > 10$, that is, $t \geq 11$. In the second case, $x > 48\frac13 - t$, forcing  $t > 19$, contradiction.

\Rf [H1] D Handelman, Positive polynomials, convex polytopes, and a random walk problem. Lecture Notes in Mathematics 1082 (1986) Springer-Verlag 133 + x.

\Rf [H2] D Handelman, Positive polynomials and product type actions of compact groups. Memoirs of the American Mathematical Society  310 (1985)  79 + xi.

%% file: generic_macros



\input colordvi
\long\def\Coloration#1#2{\textColor {#1}{#2}\textColor{.0 .0 .0 1.0}}

\def\red#1{\Coloration{.0 .99 .87 .14}#1}


\comment

olive-green  7,0,95,35       green (51,0,96,39)       red (0,100,87,14)
    darkblue (94,100,0,32)   blue (98,30,0,28)      darkbrown (0,32,61,82)
    brown (0,43,83,58)      purplered (5,98,0,29)         purple (44,97,0,44)
     verydarkblue(96,97,0,84)     pink (0,28,6,4)         orange (0,53,99,31)
\endcomment




\font\rm=cmr10 \rm

\font\bf=cmb10
\font\Rm=cmr9 at 11pt
\rm
\font\it=cmsl9 at 10pt
 at 7pt

\font\Rrm=cmr17 at 16pt
   \font\Rm=cmr12 at 11.5pt

\long\def\Pf{\par\noindent {\it Proof.} }
\def\({\left(}
\def\){\right)}
\def\st{such that }
\def\qed{\hfill$\bullet$\vskip 4pt}

\def\brcs#1{\left\{ #1\right\}}

\def\wrt{with respect to }
\def\:{\,:}

\def\aa{\text{\bf a\,}}

\def\R{\text{\bf R}}
\def\N{\text{\bf N}}
\def\Z{\text{\bf Z}}

\def\Arrow #1;#2.{#1\:#2 \to }

\def\Set#1#2{\brcs{#1 \left|\vphantom{#1 #2} \right.#2}}



\def\Rrr#1,#2{{\Cal J}_{#1,#2}}
\def\slfrac#1#2{{\raise -.07 ex\hbox{$^{#1}$}}\!/\raise .35 ex \hbox{${}_{#2}$}}
\def\ssf #1/#2{\slfrac {#1}{#2}}

\def\pd #1,#2.{\frac {\partial #1}{\partial #2}}

   \long\def\Lem
#1.#2\par{\vskip4pt{\baselineskip=13pt\font\it=cmsl12 at
11.5pt\Rm
   \noindent {\rm \uppercase{#1}} #2\vskip3pt

   }} 

\long\def\Proclaim #1.#2 \endproclaim{\vskip4pt{\baselineskip=13pt\font\it=cmsl12 at
11.5pt\Rm
   \noindent {\rm \uppercase{#1}} #2\vskip3pt

   }} 

\long\def\remark #1\endremark{\vskip 2pt \noindent {\it Remark\/} #1\par}

\long\def\Sectionhead #1.#2:\par #3{\vskip 4pt \noindent {\bf #1 #2}vskip 2pt\noindent\nospace #3}

\long\def\Title #1\par {\noindent{\Rrm #1}\vskip 9pt}

 \long\def\SubT #1.{\noindent {\it #1\/} } 
 
 \long\def\SecT
#1\par{\vskip 3pt \noindent {\bf #1}\vglue1pt
   \noindent}

\long\def\subtitle #1.{\vskip 2pt \noindent {\it #1}}

\long\def\Rmk#1\par{\vskip 1pt \noindent {\it
Remark.} #1\vskip2pt}

\long\def\Abstract #1\par{{\leftskip= 3 true cm \rightskip = 3 true cm \font\it=cmsl10 \font\rm=cmr10 \baselineskip = 10pt
\parindent=.35 true cm\rm\noindent 
{\it Abstract} #1\vskip 8pt

}}

\long\def\Author #1 \par{\noindent{\it #1}}

%% file: addonstoemptiness

\comment
Notation: $2 \leq a \leq x \leq y$; $x = ka + w$ where $0\leq w \leq a-1$; for integers $r \geq 2$ and $s$, the notation $\brm r,s. = m$   where $s\equiv m \pmod r$ and $m \in \brcs{1,2,\dots, r}$ (not the usual $\brcs{0,1,2,\dots, r-1}$).

We say an $n-1$-tuple $(a(1), \dots, a(n-1))$ (where $a(i)$ are integers strictly exceeding $1$; we also require that that they be arranged in increasing order) is {\it asymptotically hollow\/} if it satisfies the equivalent conditions of  xxx.  By that result and xxx ,

\endcomment

\SecT 4 Higher dimensional asymptotically hollow tuples with large minimum entry

As usual, we order the $n-1$-tuple $\pmb a = (a(i))$, $a(1) \leq a(2) \leq a(3) \leq \dots \leq a(n-1)$. What is the largest $a(1)$ (that is, the largest integer that can appear as the smallest integer in the sequence) that can occur in an \ah-sequence? For example, we will see later that if $n = 4$, then the largest $a(1)$ is $6$. Here we will show that for general $n$, the largest minimal entry of an \ah\ tuple (with $n-1$ entries) is {\it at least\/} $2^{n-3}(2^{n-2}-1)$. The corresponding list includes all even perfect numbers (which occur precisely when $n-2$ is a Mersenne prime, that is, when $2^{n-2}-1$ is prime).

For $n \geq 4$, define 
$$
\pmb a_n = \(2^{2n-5} - 2^{n-3}, 2^{2n-5} + 2^{n-3}, 2^{2n-4} -1, 2(2^{2n-4} -1), \dots, 2^{n-4}(2^{2n-4}-1)\) \in \Z^{n-1};
$$
 the fourth and subsequent terms (if $n \geq 5$) are double their immediate predecessor. For $n =4$, the corresponding triple is $(6,10,15)$; if $n= 5$, then the quadruple is $(28,36,63,126)$, and the next one is $(120,136,255,510,1020)$. I don't know whether these are optimal, for any $n \geq 5$.  There are an enormous number of similar constructions, but this particular one seems to maximize the least element. 

\Lem  Proposition \sixone. $\pmb a_n$ is \ah.

\Pf We use the criteria in Lemma \fivsix.  It is simple to verify that $a(1) + a(2) = a(3) + 1$, and $\sum_{i < j} a(i) = a(j)+1$ for $j \geq 3$. Moreover, we also have $\sum_{i \neq 1} a(i) = (2^{2n-4} -1)(2^{n-3}-1) + 2^{2n-5} + 2^{n-3}$. This simplifies to  
$$\eqalign{ 2^{3n-7} - 2^{2n-4} - 2^{n-3} + 1 + 2^{2n-5} + 2^{2n-3} & = 1 + 2^{3n-7} - 2^{2n-5} \cr & = 1 + (2^{2n-5} - 2^{n-3}) 2^{n-2}\cr & \equiv 1\pmod {a(1)}.\cr }$$

Finally,
$$\eqalign{
\sum_{i\neq 2} a(i) &  = \sum_{i \neq 1} a(i) + a(1) - a(2) \cr & =1 + (2^{2n-5} - 2^{n-3}) 2^{n-2} + 2\cdot 2^{2n-5} \cr  & =1 + (2^{2n-5}+ 2^{n-3}) 2^{n-2} \cr & \equiv 1 \pmod {a(2)}.\cr  }$$\qed

When $n=4$, there is only one \ah\ triple whose least entry is $6$, specifically $(6,10,15)$ (this will be shown in  section \SSS). So the following is plausible.

\Lem Conjecture \sixtwo. Suppose that $\pmb a = (a(i)) \in \Z^{n-1}$ is \ah.
\item{(i)} Then $\min_i a(i) \leq 2^{n-3}(2^{n-2}-1)$. {\par}\item{(ii)}If $\min_i a(i)  = 2^{n-3}(2^{n-2}-1)$, then $\pmb a = \pmb a_n$. 

The first test of this conjecture occurs when $n =5$; for (i), we would have to show that $a(1) \leq 28$. And the natural  conjecture concerning the largest $a(i)$:  $a(n-1)$ should be at most the largest coefficient of $\pmb a_n$.

\SecT 5 Proscriptive intervals

One approach to the problem of characterization of \ah\ sequences yields intervals of real numbers in which the last coordinate cannot appear. Specifically, let $(a(1), a(2), \dots , a(n-2))$ be   monotone increasing. We wish to determine (among other things) which integers $ y \geq a(n-2)$ (or whether such a $y$ exists) \st $\pmb a := (a(1), a(2), \dots, a(n-2), y)$ is asymptotically hollow.   We obtain fairly good constraints on $y$. As a preliminary example, suppose that $ y \geq a(n-2)t$ for some positive integer $t$. Then for all $i \leq n-2$, we have $a(i)t \leq a(n-2)t \leq y$, and thus $\brm y,a(i)t. = a(i)t$. Hence $(n-3)y + t \geq \sum_{i=1}^{n-2} a(i)t$. Set $s = \sum_{i=1}^{n-2}a(i) -1$. 

Therefore, $y \geq a(n-2)t$ implies $ y \geq ts/(n-3)$. That is, $y$ does not belong to the half-open real intervals 
$\left[a(n-2)t, st/(n-3)\)$ for all positive integers $t$. Hence
$$
y \notin \bigcup_{t \in \N} \left[a(n-2)t, \frac{st}{n-3}\) =  \bigcup_{t \in \N} t\left[a(n-2), \frac{s}{n-3}\).
$$
It can happen   that $a(n-2) \geq s/(n-3)$, in which case all of these multiple intervals  are empty, and no constraint on $y$ results . On the other hand, if $a(n-2) < s/(n-3)$, then the union contains an infinite ray, that is, there exists $b \in \R$, \st $[b, \infty)$ is contained in the union. It is easy to estimate $b$, for example, determine the smallest integer $t_0$ \st the right endpoint of $st_0/(n-3) \geq a(n-2)(t_0 +1)$, and then we can set $ b= a(n-2)t_0$ (smaller choices for $b$ can sometimes arise, especially since we only care about integer values for $y$). This yields a first class of constraints. There are  more of this type. 

To explain this, order the set $\Cal S:= \Set{a(i)/m}{ i \leq n-2; m \in \N}$ in the usual way (as a set of real numbers; generally, ties occur, that is, $a(i)/m = a(i')/m'$ can arise with differing $i,i'$). For each $a(i)/m$, define $g(i,j.m) = \#\Set{k \in \N}{a(i)/m < a(j)/k}$. If $i = j$, then $g(i,i,m) = m-1$, and obviously $g(i,i,1) = 0$. In general, $g(i,j,m)$ is simply the largest integer $k$ \st $k < ma(j)/a(i)$, that is,  $$\eqalign{ g(i,j,m) &= \flo{\frac{ma(j)-1}{a(i)}} \cr  &= \frac{ma(j) - \brm a(i),ma(j).}{a(i)}. \cr }$$

Define $f(i,m) = \sum_{j=1}^{n-2} g(i,j,m)$. Now let $c \equiv c(i,m) \in \Cal S$ be the smallest \st $ c> a(i)/m$ (if no such $c$ exists, then $m=1$ and $a(i) = a(n-2)$, and the resulting proscriptive intervals will be the ones obtained above from $ y \geq a(n-1)t$, that is, nothing new happens; so we may assume $c$ exists). 

Suppose that $a(i)t/m \leq y < ct$ for some $t \in \N$. Then $a(i)t/m \leq y < a(i)t/(m-1)$ (true even if $ m=1$), so that $(m-1)ty < a(i)t \leq mty$. Thus $\brm y,a(i)t. = a(i)t - (m-1)y$, and this is $a(i)t - g(i,i,m)y$. Now let $ j\neq i$. If $a(j) \leq a(i)/m$, then $a(j)t \leq y$, and $\brm y,ta(j). = ta(j)$, and of course, this is $ta(j) -g(i,j,m)y$. On the other hand, if $a(j) > a(i)/m$, then $a(j) \geq c$, and thus $g(i,j,m) \geq 1$. 

We have that $a(j)/g(i,j,m) \geq a(i)/m$ (by definition), with equality only if $a(j)/g(i,j,m) = a(i)/m$, and thus $\brm y,a(j)t . = a(j)t - g(i,j,m) y$ in this case as well. Finally, if $a(j)/g(i,j,m) > a(i)/m$, then $a(j)/g(i,j,m) \geq c$, and thus $a(j)t/g(i,j,m) \geq ct > y$. This yields $a(j)t > g(i,j,m)y$. 

On the other hand, $a(j) < (g(i,j,m) +1)a(i)/m$, so that $a(j)t <(g(i,j,m) +1) a(i)t/m \leq (g(i,j,m)+1)y$. Hence $\brm y,ta(j) . = ta(j) - g(i,j,m)y$ in all cases. Thus  for all $i$,
$$\eqalign{
\sum_{j=1}^{n-2} \brm y,ta(j). & =\sum_{j=1}^{n-2} ta(j) -  y\sum_{j=1}^{n-2} g(i,j,m)\cr & = (s+1)t - f(i,m)y \cr 
}$$
From  $  \sum_{j=1}^{n-1} \brm y,ta(j). \leq (n-3) y + t$, we deduce  
$$
 y \geq \frac{st}{n-3 + f(i,m)}. 
$$
Set $H(i,m) = n-3 + f(i,m)$. We obtain the corresponding {\it proscriptive intervals,}\plainfootnote{$^1$}{{\it Pr\red{{o}}scribe} is the opposite of {\it pr\red{{e}}scribe}.}%
 (sometimes, {\it proscribed intervals})  
$$
 y \notin \left[\frac{ta(i)}{m}, \frac{st}{H(i,m)} \wedge c(i,m)t \right). 
$$ 
The infimum, $\wedge$, arises from the original condition that $y < ct$, but we can dispense with the $c(i,m)$ term, as follows. If $c't > y\geq c(i,m)t$ (where $c'$ is the next strictly larger term than $c$ (ignore ties) in $\Cal S$), then we similarly derive, $ y\not\in [ct,st/H(i_0,m_0)\wedge c't)$, where $ c = a(i_0)/m_0$. Then from the definitions, $H(i,m) < H(i_0,m_0)$, and we continue by induction on the elements of $\Cal S$ increasing to $a(n-1)$, to eliminate the wedge term. In particular, the corresponding proscriptive interval is $\left[\frac{ta(i)}{m}, \frac{st}{H(i,m)}   \right)$, and we have 
$$
y \notin \bigcup_{t\in \N}\left[\frac{ta(i)}{m}, \frac{st}{H(i,m)} \right). 
$$
The proscriptive intervals are determined  by $t$, $a(i)$,   $m$, and $s$, and we allow $m$  and $t$ to vary over all positive integers. 

We do not require $y \geq a(n-2)$ in this construction, and this will be useful in problems of the type discussed in Example \sevfou. 

There are only finitely many choices of $(i,m)$ for which the interval $[a(i)/m, s/H(i,m))$ is nonempty. To see this, expand 
$$\eqalign{ H(i,m)& = n-3 + m-1 + \sum_{j\neq i} g(i,j,m)\cr  
&= n- 3 + m-1 + \frac{m}{a(i)} \sum_{j\neq i} a(j) - \frac 1{a(i)}\sum_{j\neq i} \brm y,ma(j).\cr 
&= n+ m-4 + \frac{m}{a(i)} (s - a(i)+ 1) - \frac 1{a(i)}\sum_{j\neq i} \brm y,ma(j).\cr 
&= n-4 + \frac{m(s+1)}{a(i)} - \frac 1{a(i)}\sum_{j\neq i} \brm y,ma(j). .\cr  
}$$ 
The interval $[a(i)/m, s/H(i,m))$ (and thus all of its  positive integer multiples) is empty if and only if $a(i)H(i,m)/m \geq s$. Applying the preceding expression for $H(i,m)$, the interval is trivial if and only if  
$$ \eqalign{ s+1 + \frac{(n-4)a(i)}{m} &- \frac{1}{m}\sum_{j\neq i} \brm y,ma(j). \geq s, \quad\text{that is}\cr  \sum_{j\neq i} \brm y,ma(j). &\leq m + (n-4)a(i)\cr 
}$$
 These are  stronger than the necessary conditions (for $\pmb a$ to be asymptotically hollow), that $\sum_{j\neq i} \brm y,ma(j). \leq m + (n-3)a(i)$, and lends itself immediately to the following. 

\Lem Corollary \sevone. If all proscriptive intervals of $\pmb b:= (a(1), a(2), \dots, a(n-2))$ are empty, then $\pmb b$ is asymptotically hollow. 

\Pf Simply note that then $n-2$-tuple corresponds to $n$ replaced by $n-1$, and $(n-1)-3 = n-4$.

\Lem Corollary \sevtwo. If there exist infinitely many $y$ \st $(a(1), a(2), \dots, a(n-2),y)$ is asymptotically hollow, then $(a(1), \dots, a(n-2))$ is asymptotically hollow. 

\Pf If at least one of the proscriptive intervals is nontrivial, then its union over all positive integer multiples contains an infinite ray, and thus there is an upper bound on the possible $y$ \st $(a(1), a(2), \dots, a(n-2),y)$ is asymptotically hollow. \qed 

More can be deduced from the methods here. Fix $i,j$  and let $k = g(i,j,m)$. Then $k+1 \geq ma(j)/a(i) > k$ by definition, and thus $(k+1)a(i) \geq ma(j) > ka(i)$. Thus $\brm a(i),ma(j). = ma(j) - g(i,j,m) a(i)$. Allowing $j$ to vary, we deduce 
$$\eqalign{
\sum_{j \neq i}\brm a(i),ma(j). & = m \sum_{i\neq j} a(j) - \sum_{i\neq j} g(i,j,m) a(i)\cr & = m(s+1-a(i)) - \(f(i,m) - (m-1)\)a(i) \cr & = ms - f(i,m) + 1)a(i); \quad \text{thus}\cr  \brm a(i),my. &\leq (n-3)a(i) + m - \sum_{j\neq i} \brm a(i),ma(j). \cr  & = (n-2 + f(i,m))a(i) - ms  = (1 + H(i,m))a(i) - ms. \cr 
}$$ 
This imposes constraints on $\brm a(i), y.$, which are sometimes useful (more frequently with $i = n-2$). The condition for triviality of the corresponding interval, $[a(i)/m, s/H(i,m))$ is that $ms \leq H(i,m) a(i)$, which is slightly stronger than the inequality here. 

We also have a quantitative proof of the preceding, which relates proscriptive intervals to solutions to inequalities.

\Lem Proposition \sevthr. Suppose that $\pmb a':= (a(1), a(2), \dots, a(n-2)) \in \N^{n-2}$ is not \ah.  \item{(a)} There exist positive integers $i$ and $m$ with $i \leq n-2$ and $m < a(i)/2$ \st $$
\sum_{l\neq i} \brm a(i), a(l)m. > (n-4) a(i) +m ,\tag1 $$ and the corresponding proscriptive interval $[a(i)/m, s/H(i,m))$ (where $s = \sum_{j=1}^{n-2} a(i) - 1$) is nonempty. 
 \item{(b)} Suppose that $i,m$ are positive integers \st (1) holds.  If $\pmb a:= (a(1), \dots, a(n-2),y) \in \N^{n-1}$ is \ah\ for some positive integer $y$, then  
$$ 
y  < \frac{a(i)(1 + s-m)}{2m}. 
$$

\Pf The first part of (a) is a consequence of the general criterion for asymptotic hollowness, applied with $n-1$ replacing $n-2$. Consequent nonemptiness of the interval was shown earlier.  \noindent (b) Set $I = [a(i)/m, s/H(i,m))$, so that $y \not\in \cup_{t\in \N}tI$. There exists $t_0 \in \N$ \st for all $t \geq t_0$, we have $st/H(i,m) \geq a(i)/m$, estimated as follows. The inequality can be rewritten 
$$\eqalign{
t_0\(\frac{ms}{H(i,m)} - a(i) \)& \geq a(i); \quad{\text{so}}\cr 
t_0 & \geq \frac{a(i)H(i,m)}{ms - a(i)H(i,m)} \cr 
&= \frac{a(i)(n+m-4) + m\sum_{l\neq i} a(l) - \sum_{l\neq i}\brm a(i),ma(l). }{ms - \(a(i)(n+m-4) + m\sum_{l\neq i} a(l) - \sum_{l\neq i}\brm a(i),ma(l).\)}\cr
 &=  \frac{a(i)(n-4) + m + ms  - \sum_{l\neq i}\brm a(i),ma(l). }{m- a(i)(n-4)  +\sum_{l\neq i}\brm a(i),ma(l).}\cr  
&= \frac
{\frac{a(i)(n-4)}m + s+1  - \frac{\sum_{l\neq i}\brm a(i),ma(l).}m }
{1- \frac{a(i)(n-4)}m   +\frac{\sum_{l\neq i}\brm a(i),ma(l).}m}.\cr 
}$$

If $t_0$ is set to the ceiling of the last expression, than any larger $t$ will also satisfy the inequality, and thus the ray $[a(i)t_0/m, \infty)$ is contained in $\cup_{t\in \N} tI$. In particular, $y \leq a(i)t_0/m$. The denominator is at least as large as $2$ (by (1)), and the numerator is bounded above by $s+1 - m$. So we can take $t_0 = \ceil{(s+1 - m)/2}$, e.g., $t_0 \leq (s-m)/2 + 1$. Thus $a(i)t_0/m \leq (1 + (s-m))a(i)/2m$.  \qed \par \par If there are numerous pairs $(i,m)$ \st (1) holds, then the the upper bound can be made much smaller (often), since there will be overlapping of multiples of different intervals.

On rereading this article for the $(n-1)$st time, I decided that the definition of $H(i,m)$ hides its essential simplicity. We provide an example to show how to use it effectively. 

 \Lem Example \sevfou. Suppose $\aa$ is a nontrivial \ah\ quadruple that has entries $29, 38, 66$. Then $\aa$ is one of  $(2, 29, 38, 66)$, $(3, 29, 38, 66)$, or $(11, 29, 38, 66)$. 

 \Pf We  consider the quadruple $(29, 38, 66, y)$, where we do not require $y \geq 66$. Here $s = 132$. We obtain various proscriptive intervals, out of the saturated  sequence of fractions 
$$
 \frac {66}1 > \frac{38}1 > \frac{66}2 > \frac{29}1 > \frac{66}3 > \frac{38}2 > \frac{66}4 > \frac{29}2 >\frac{66}5 > \frac{38}3 > \dots.
$$
 This initial sequence yields the respective proscriptive intervals (we will need more later), 
$$\eqalign{
 t\left[66, \frac{132}2\) = \emptyset; \quad t\left[38, \frac{132}3\) = t\left[38,44\); 
& \quad t\left[33,\frac{132}4\) = \emptyset; \quad t\left[29,\frac{132}5\) = \emptyset; \cr
\quad t\left[22,\frac{132}6\) = \emptyset; \quad t\left[19,\frac{132}{7} \)= t\left[19, 19\frac{5}7\);& \quad t\left[15\frac{1}2, \frac{132}8\) = \emptyset;\qquad \text{etc}. 
}$$
 The denominators, $H(i,m)$, just increase by one.  
 
 Now consider what happens modulo $38$: we have $\brm 38,29t. + \brm 38,28t. + \brm 38,ty. \leq t + 77$ (as $66 \equiv 28 \pmod {38}$) for all $t$, but just from $t = 1$, we have $\brm 39,y. \leq 79-67 = 11$. In particular, $ y = 38k + j$ where $k \geq 0$ and $1 \leq j \leq 11$ are integers. From the proscriptive intervals $[38,44) \cup [76,88) \cup [114 , 132) \cup \dots$, we see that the only choices for $y$ are $\brcs{2, \dots, 11} \cup \brcs{45,\dots, 49}$. 
 
 Now consider the situation modulo $29$; we have $\brm 29,9t. +\brm 29, 8t. + \brm 29,ty. \leq t + 58$. With $t = 3$, we have $27 + 24 + \brm 29,3y. \leq 61$, and so $\brm 29,3y. \leq   10$. The only possible values (of the previous paragraph) satisfying this are $ y \in \brcs{2,3,10,11, 49}$. We can eliminate $10$ and $49$ simultaneously. 
 
 The saturated sequence of fractions around $10$ is 
$$ \frac{29}2 > \frac{66}5 >\frac{38}3 >\frac{66}{6} > \frac {29}3 > \frac{38}4 > \dots.
$$
 To find $H(1,3)$  (here $a(1) = 29$), we simply add the denominators for each of the terms in $29$ $(3)$, $38$ (3), $66$ (6), and add one (because $ n =5$), so that $H(1,3) = 13$ (the occurrence of $66/5$ between $38/3$ and $29/2$, but before $66/6$, is irrelevant for this). Thus we obtain the proscriptive intervals 
$$t\left[\frac{29}3, \frac{132}{13}\) = t\left[9\frac 23, 10\frac2{13}\).
$$
 At $t=1$, we obtain that $10$ belongs to the interval, so that $ y \neq 10$, and at $t = 5$, $49$ belongs, eliminating it. We are thus reduced to $ y \in \brcs{2,3,11}$. 

 [As an aside, we also see directly, that if $t = 5$, then modulo $49$ (for the quadruple $(29,38,49,66)$) we have $\brm 49,29\cdot 5. + \brm 49,38\cdot 5. + \brm 49,66\cdot 5. > 5 + 98$, either by direct computation or from Lemma \sevthr.] 
 
 Now we show that with each of $y = 2, 3, $ or $11$, the resulting quadruples are \ah. We notice that $38 + 29 \equiv 1 \pmod {66}$, so we need only consider the inequalities modulo $y$, $29$, or $39$. 

 \vskip2pt\noindent $\bullet$ {\it $(2, 29, 38, 66)$ is \ah}. Modulo $2$ is trivial. We can deal with both $29$ and $38$ because of the following easy result.  

 \Lem Lemma \sevfiv. Let $(2,b,c,d)$ be a quadruple of integers exceeding $1$, not necessarily in ascending order. Suppose that $c \equiv b + 1 \pmod d$. Let $t$ be a positive integer less than $d/2$. 
Then  $\brm d,2t. + \brm d, bt. + \brm d, ct. \leq t + 2d$.

 \Pf Assume that $\brm d,2t. + \brm d, bt. + \brm d, ct. > t + 2d$. Set $ x = \brm d,bt.$. Then $\brm d,ct. = \brm d,(b+1) t. \equiv x + t \pmod d$. Also $2t \leq d$, so $\brm d,2t. = 2t$. 

 If $x+t \leq d$, then $\brm d,ct. = x+t$. Then we have 
$$
 \brm d,2t. + \brm d, bt. + \brm d, ct. =2t + 2x + t .
$$
 Plugging these into the inequality, we obtain 
$$
  2(x+t) > 2d,  
$$
contradicting $x + t \leq d$. 

 Now suppose that $x+ t > d$; as $x \leq d$ and $t < d$, it follows that $x+t < 2d$, so that $\brm d,ct. = x+t -d$. This yields $2t + 2 x + t - d > t + 2d$, so $2x > 3d - 2t$. As $2x \leq 2d$, we have $d < \brm d,2t.$, a contradiction. \qed 

 This applies directly to $y =2$ and both $d =29$ and $38$, since modulo $29$, the remainders are $2,9,8$, and modulo $38$, the remainders are $2,29,28$. This shows that $(2,29,38, 66)$ is \ah.\qed 

 \noindent $\bullet$ {\it $(3,29,38,66)$ is \ah.}  First we consider the case that $d = 29$, that is the inequalities are verified modulo $29$. We reduce modulo $29$ to $(3, 8, 9)$---we wish to verify that for all $1 \leq t \leq 14$, we have $\brm 29,3t. +  \brm 29,8t. + \brm 29,9t. \leq t + 58$. This was proved just after Lemma \fivsev, but we want to illustrate the use of proscriptive intervals.  Assume that for some  $t$,  $\brm 29,3t. +  \brm 29,8t. + \brm 29,9t. > t + 58$. 

 Denote   $\brm 29,8t.$ by $x$, so that $\brm 29,9t. \equiv x+ t \pmod {29}$. Since $29$ is prime and $t \leq 14< 29$, we must have $x \leq 28$. We note that $3 + 8 + 2\cdot 9 = 29$, and so $29$ divides $3t + x + 2x+2t = 3x + 5t$. First suppose that $t\leq 9$, so that $\brm 29,3t. = 3t$. If $x + t > 29$, then (as $x \leq 28$), $\brm 29,9t. \leq t-1$. Then we obtain $3t + x + t-1 > 58 +t$, so $x + 3t > 59$; as $t \leq 9$, it follows that $x > 32$ a contradiction. Hence $x + t \leq 29$, and thus $\brm 29,9t . =x+t \leq 28$. This yields $3t + x + x+ t > t+ 58$, so $2x + 3t > 58$. This entails $3x + 9t/2 > 87$, and thus $3x + 5t > 87$; the left side being divisible by $29$ forces $3x + 5t = 29\cdot k$ where $k \geq 4$ is an integer. On the other hand, $3x +5t \leq 84 + 45 = 129 < 145$. Hence $3x + 5t = 116$. Thus $5t \geq 116 - 84 = 32$, so $7 \leq t \leq 9$. Also, $3x +5t =116$ entails $t \equiv 1 \pmod 3$, and thus $t = 7$, forcing $x = 27$, which contradicts $x + t \leq 29$. 
 
 Now suppose that $10 \leq t \leq 14$; then $\brm 29,3t. = 3t - 29$. If $x + t > 29$, as before it follows that $\brm 29,9t. \leq t-1$, and so $3t - 29 + x + t-1 > 58 +t$, that is, $x + 3t > 88$. Since $3t \leq 42$, we obtain an immediate contradiction. Thus $x + t \leq 29$, and so $3t-29 + x + x+t > 58 + t$, that is, $2x + 3t > 87$. As before, $3x + 9t/2 \geq 132$, so that $3x + 5t = k\cdot 29$ where $k \geq 5$ is an integer. As $3x + 5t \leq 84 + 70 = 154$, we have that $3x + 5t = 145$. Since $2x + 3t \leq 87$, we have $x + 2t \geq 58$, contradicting $x + 2t \leq 28 + 28$.   
 
 Now we consider $(3,29,38, 66)$ modulo $38$. This reduces to $(3, 28, 29)$. Assume $\brm 38,3t. + \brm 38,28t. + \brm 38,29t. > t + 76$ for some positive integer $t \leq 18$, and set $x = \brm 38,28t.$. Since $2\cdot 29 + 2 \cdot 28 = 114$, we have  $38$ divides $4x + 2t$. 
 
 First, consider the case $t \leq 12$. Then $\brm 38,3t. = 3t$. If $x + t > 38$, then $\brm 38,29t. \leq t$, and so $3t + x + t > t + 76$, that is, $x +3t > 76$, yielding a contradiction, as $3t \leq 36$. Hence $x + t \leq 38$, and thus $\brm 39,29t. = x+t$. 

 Then $3t + x + x+t > t + 76$, that is, $2x + 3t > 76$, so $4x + 6t > 152$; also, $4x + 6t \leq 4\cdot 38 + 6\cdot 18  <7\cdot 38$. Write $4x + 2t = 38\cdot k$ for some $k$. Then $7\cdot 38 >4t + 38 k > 152$; the right side yields $k \geq 3$. If $ k =3$, then $4t > 38$, forcing $10 \leq t \leq 12$. Since $2x + t = 57$, $t$ is odd, yielding $t= 11$. But then $x = 23$. However, $x \equiv 28\cdot 11 \equiv -10\cdot 11 \pmod 38 $, and this is just $3$, a contradiction.  

 If $k = 4$, then $x = 38-t/2$, contradicting $x + t \leq 38$. If $k \geq 5$, then $x \geq 95/2 - t/2 \geq 83/2 > 38$, contradiction. 
 
 So we are left with the case that $13 \leq t \leq 18$, and then $\brm 38,3t. = 3t - 38$. If $x + t > 38$, then as before, $\brm 38,29t. \leq t$ and so $3t - 38 + x + t > t + 76$, that is, $x + 3t > 114$, whence $x > 60$ a contradiction. Thus $x + t \leq 38$. 

 This yields $3t-38 + 2x + t > t + 76$, so $2x + 3t > 114$. But $2x + 2t \leq 76$ and $t \leq 18$, a contradiction. \qed  

 \noindent $\bullet$ {\it $(11, 29, 38, 66)$ is \ah}. We have $29 + 38 \equiv 1\pmod {11}$ and $11 + 66 \equiv 1 \pmod {38} $; by Lemma \fivsix, we need only deal with one situation, that is,  modulo $29$. 
 
 As usual, assume that $\brm 29,11t. + \brm 29,8t. + \brm 29,9t. > t + 58$ for some positive integer $t \leq 14$, and set $x = \brm 29,8t.$. We observe that $\brm 29,9t. \equiv x + t \pmod {29}$ and $\brm 29,11t. \equiv x + 3t \pmod {29}$.  In addition, $11 + 38 + 66= 115 \equiv -1 \pmod {29}$, and so $x + 3 t + x + x+t \equiv -t \pmod{29}$; alternatively, $3x + 5t$ is divisible by $29$. Now $3x + 5t < 3\cdot 29 + 70 = 147$ (and less than $145$ if $t \neq 14$. Hence, on writing $3x + 5t = 29k$, we must have $k \leq 6$, and if $t \neq 14$, then $k \leq 5$. 
 
 If $x + t > 29$, then $\brm 29,9t. \leq t-1$, and we deduce (over-estimating $\brm 29,11t. $ by $28$), that $28 + x +t-1 > t + 57$, that is, $x > 30$ a contradiction. Hence $x + t \leq 29$, so $\brm 29,9t. = x+t$, and this is at most $28$ (since $29$ is prime and $t \leq 14$). If now $x+ 3t > 29$, then $\brm  29,11t. \leq 2t-1$. This yields $2t-1 + x + x+t > t+ 58$, so $2(x+t) > 59$, a contradiction. Therefore, $ x + 3t \leq 29$, so $\brm 29,11t. = x+3t$, and again by primeness of $29$, $x+ 3t \leq28$. 

 We obtain $x + 3t + x+ x+ t > t +58$, or $3x + 3t > 58$. As $x + 3t \leq 28$, $3x + 9t \leq 84$; hence $3x + 5t < 84$, and therefore $k\leq  2$, that is, $3x + 5t \leq 58$, which obviously contradicts $3x + 3t > 58$. \qed

\SecT \SSS \ Classification when $ n=4$

In this section, we will analyze the possibilities for $n = 4$. It turns out that there is one infinite family of (nontrivial) asymptotically hollow triples, and a small number of others. We relabel the triples, $(a,x,y)$  with $a \leq x \leq y$ (instead of $(a(1), a(2),a(3))$) for convenience. Then the conditions of Corollary \fivfiv\ translate to the following three inequalities,   
$$\eqalign{ \brm a,xt. + \brm a,yt.&\leq t + a \quad \text{for all positive integers $t < a/2$}\cr 
\brm x,at. + \brm x,yt.&\leq t + x \quad \text{for all positive integers $t < x/2$}\cr 
\brm y,at. + \brm y,xt.&\leq t + y \quad \text{for all positive integers $t < y/2$.}\cr }$$

Suppose that $(a,x,y)$ is asymptotically hollow, and $a \geq 2$. Define  $k = \flo{x/a}$; then  
$x= ka + \delta$, where $k \geq 1$ and $0 \leq 
\delta \leq a-1$. Then $x > x/2 > \dots > x/k \geq a > x/(k+1)$. We recall the construction of proscriptive intervals from the previous section.

If $y \geq xt$ for some positive integer  $t $ (necessarily with $t < y/2$), then $\brm y,xt. = xt $, and as $a \leq x$, we also deduce $\brm y,at. = at$. The third inequality yields $at + xt \leq y + t$, or in other words, $y \geq (x + a-1)t$. So $y  \geq xt$ entails that $y \geq (x+a-1)t$. Hence $y \notin  S_1:=\cup_{t \in \N}[xt, (x+a-1)t)$ (the first and easiest to deal with union of proscriptive intervals). We similarly  deal with the possibility that $xt/r \leq y < xt/(r-1)$ for some $r \in \brcs{2,3,\dots,k}$ (if $k \geq 2$). It follows that $xt - (r-1)y > 0$, but $xt - ry < 0$. Hence $\brm y,xt. = xt - (r-1)y$. Also, as $a \leq x/k$, $at \leq xt/k \leq y$, we have $\brm y,at. = at$. Applying the third inequality as before, we obtain $y \geq (x+a-1)t/r$. Since $(x+a-1)/r \leq xt/(r-1)$ (as $ka \leq x$), we have that $y \not\in S_r:=\cup_{t \in \N}[xt/r, (x+a-1)t/r)$. From $y \notin S_1$ and with $t =1$, we observe (trivially) that $y  \geq x+a-1$. 

 Now we can give  a complete description  of   asymptotically  hollow triples ($n=4$).

 \Lem Theorem \Sone. The nontrivial asymptotically hollow triples are precisely{
}\item{(i)}  $(2,x,x+1)$ for every $x \geq 2$ 
 \item{(ii)} $(2,3,5), (2,3,8),    (2,5,9)$ 
 \item{(iii)} $(3,4,6), (3,5,7), (3,5,8), (3,8,10)$ 
 \item{(iv)} $(4,6,9), (4,7,10)$ 
\item{(v)} $(5,8,12)$ 
\item{(vi)} $(6,10,15)$.  
 
 (Observe that  $(2,3,4)$ and $(2,5,6)$ are asymptotically hollow; they appear in (i).) It is an easy consequence that all the hollow lattice polytopes that arise from these asymptotically hollow triples have width at most two (these are of course four-simplices; finitely many empty four-simplices with width exceeding two exist [IS], but none of them  appear here). 
The first and easy step is to verify that all the listed triples are
asymptotically hollow. For most of the computations, we can apply Lemma\, \fivsix, which specializes to the following when $n =4$.  
 
 \Lem Lemma \Stwo. Suppose $a,b,c$ are  integers all of  which are at least two, and  any  of
the following conditions hold. \item{(a)} either $a \equiv 1 \pmod c$ or $b
\equiv 1 \pmod c$; 
\item{(b)} $a+ b \equiv 1 \pmod c$; 
 \item{(c)} $a+ b \equiv 0 \pmod c$  and at least one of $a$ or $b$ does not divide $c$. {

} \noindent Then for all integers $t$ with $1 \leq t < c/2$, we have 
$$
 \brm c,ta. + \brm c,tb. \leq t + c.
$$

For example, to verify that  $(2,x,x+1)$ is asymptotically hollow for any integer $x \geq 2$, we
observe that at least one of $x$ or $x+1$ is odd, so the condition holds for $c
=2$; next $x +1 \equiv 1 \pmod x$, and $2 + x \equiv 1 \pmod {x+1}$, so the
conditions hold for $c = x$ and $c = x+1$. Each of those whose smallest entry is at least $4$ are immediately verified to be asymptotically hollow, by three applications of Lemma \Stwo. For some of these triples, for example,
$(2,3,8)$, we have to check   some of the conditions by hand, but this is routine. Or we could just use Lemma \fivsev. 

Now we begin the proof that that the listed triples are the only nontrivial
asymptotically hollow ones. 
 We write the candidate as $(a,x,y)$, and    set $k = \flo{x/a}$; then  
$x= ka + \delta$, where $k \geq 1$ and $0 \leq \delta \leq a-1$. Also define
$\rho = y-x$. 
 
 The proof is divided into five parts (which further subdivide), roughly
speaking, as follows (more detail is included in the precise statements). 

 \item{(A)} No asymptotically hollow triples exist with $k \geq a \geq 3$. 
 \item{(B)} No asymptotically hollow triples exist with $3 \leq k < a$. 
 \item {(C)} The only  asymptotically hollow triple  with $k =2$ and $a > 2$ is
$(3,8,10)$.
 \item{(D)} The only asymptotically hollow triples  with $k \geq 2$  and $a =2$ are $(2,5,9)$ and $(2,x,x+1)$ with $x \geq 4$. 
 \item{(E)} $k=1$.

We   severely restrict the possibilities for $y$ (and consequently for $a$ and $x$) via proscriptive intervals. For $r \leq k$, let $S_r$ be the union of the proscriptive intervals 
$$
S_r:= \bigcup_{t = r, r+1, r+2, \dots} \left[\frac{xt}r, \frac{(x+a-1)t}r \). 
$$
In $S_r$, consider the condition that the right endpoint for some $t = t_0$, $(x+a-1)t/r$) is at least as large 
as the left endpoint for the prescriptive interval arising from $t_0+1$ ($x(t+1)/r))$, that is, $(x+a-1)t/r \geq x(t+1)/r$; this is simply that $t_0 \geq x/(a-1)$ and this is independently of $r$. Set $t_0 = \ceil{x/(a-1)}$. We deduce that $S_r$ contains the infinite ray, $[x t_0/r, \infty)$. The largest of these rays arises when $r = k$, and this forces  $ y < xt_0/k$. 

As a rough (for now) estimate, $t_0 = \ceil {x/(a-1)} = k + \ceil{(\delta+k)/(a-1)}$; 
hence $y \leq x (1 + \ceil{(\delta +k)/(a-1)}/k )$. We see that $\ceil{(\delta+k)/(a-1)}/k < 2/k + 1/(a-1)$. In particular if $2/k + 1/(a-1) \leq 1$, then $ y \leq 2x$; since $2x$ belongs to a proscriptive interval, in this case $ y \leq 2x-1$.

For example, if $ k \geq 3$ and $a \geq 4$, or if $ k \geq 4$ and $a \geq 3$, then $y \leq 2x$. When $k = 2$, 
we don't use this estimate, but go back to $\ceil{\delta+ 2)/(a-1)} \in \brcs{1,2}$, so obtain $y \leq 3x/2$  when $a \geq 3$, and since $3x/2$ is the initial endpoint of a proscriptive interval, we have $y < 3x/2$;  similarly if $ k = 3$, $\ceil {(\delta+3)/(a-1)} \in \brcs{1,2}$  when $a \geq 3$, so that $y-x \leq 2x/3$. If $k = 4$, we similarly obtain that if $a \geq 3$, then $\rho \leq 3x/4$.

For larger $k$ (but still with $a \geq 3$), we require in what follows the simple inequality $(y-x)/x \leq 1 - 1/k$. But this is a consequence of $2/k + 1/(a-1) \leq 1 - 1/k$, that is, $3/k + 1/(a-1) \leq 1$; this obviously holds if $k \geq 6$. When $k = 5$ and $a = 4$, $\ceil{(k+\delta)/(a-1)}/k = \ceil{(5+\delta)/3)}/5 \leq 3/5 < 4/5$; when $k = 5$ and $a = 3$, we obtain $\ceil{(4+\delta)/2}/4 \leq 3/4$. We conclude that provided $a \geq 3$ (and our running assumption that $k \geq 2$), $(y-x)/x \leq 1-1/k$. We will further restrict the the possible values of $(y-x)/x$ as we go along, by considering other $S_r$  (when $k \geq 2$). We summarize these remarks.

\Lem Lemma \Sthi. (a) $y \leq x(1 + 2/k + 1/(a-1))$. 
\item{(b)} If $k \geq 2$ and $a \geq 3$, then $y < (2-1/k)x$.
\item{(c)} $y < \ceil{x/(a-1)}\cdot x/k$.

 \Lem (A) Lemma \Sfou. No \ah\ triples exist with $k \geq a \geq 3 $. 

We require a preliminary lemma. 
 
 \Lem Lemma \Sfiv. If $a| x$, then $a =2$,  and in this case, the only asymptotically hollow triples are of the form $(2,2m,2m+1)$  ($m=1,2,\dots$). 
 
 \Pf \noindent     Here the triple is $(a,ka,y:= ak + \rho)$; $\delta = 0$, $t_0 = \ceil{ak/(a-1)} = k + \ceil{k/(a-1)}$, and $y < [x(1 + 2/k + 1/(a-1)), \infty)$; the left endpoint is  $x + 2a + ka/(a-1)$. Also $ y \notin [at,(x+a-1)t/(k+1)]$. Moreover, modulo $a$, we deduce $a +\brm a, \rho. \leq a+1$, so $\brm a,\rho. = 1$. First, suppose    $k = 1$, so that $t_0 = 2$, and $ y \not\in [2x, \infty)$, that is, $ \rho < x$. Then applying the criteria modulo $x$, we have $a + \rho \leq x+1 = a+2$, so that $\rho \leq 2$. The resulting triples are thus $(a,a,a+1)$ or $(a,a,a+2)$. In both cases, modulo $a+1$, we have $a + a \leq a+2 \text{ or } a+3$, respectively. The former yields  $a = 2$, the latter $a \leq 3$. With $a= 2$, we have the triple $\boxed{(2,2,3)}$, which despite the repetition, is asymptotically hollow; the other triple is $(2,2,4)$ which is obviously not asymptotically hollow. The remaining case is $(3,3,5)$, which violates $\rho \equiv 1 \mod 3$.

Now let $k \geq 3$. We have $ y \notin [xt/k, (x+a-1)t/k)  = [at,at + (a-1)t/k)$. If $a \geq 4$, then $1/2k + 1/(a-1) \leq 1$, and so $\rho < x$; thus $\rho = la + 1$ with $l < k$. Set $t = l+k$; then $y \not\in [a(l+k), a(l+k) + (a-1)(l+k)/k$. Since $(a-1)(l+k)/k > 1$, we are done. If $k \geq 3$ and $a = 3$, then $x = 3k$, and $t_0 = \ceil{x/2} \leq (x+1)/2$, so $ y \notin [3(3k+2)/2, \infty)$, so $\rho < (3k+2)/2 $; this is less than $x = 3k$. Write $\rho = 3l+1$; necessarily, $l < k$. Set $t = l+k$ as before; we obtain $y \notin [3(l+k), 3(l+k) + 2(l+k)/k)$, and since $2(l+k)/k > 1$, we are done.

When $a = 2$, we have $(2,2k,2k+\rho)$, and now $y \notin [2t, 2t + t/k)$; here $t_0 = x$, so $ y \not\in [x^2/k, \infty) = [4k,\infty) = [2x,\infty]$. Thus $\rho < x$, and we can again set $t = l+k$ where $\rho = 2l+1$ (with $ l < k$). If $l \neq 0$, we obtain $ y \notin [2l, 2l + (l+k)/k)$, and once again  $(l+k )/k > 1$, ruling out all possible $\rho$.  The only remaining cases are of the form  $\brcs{(2,2k,2k+1)}$, which are asymptotically hollow. Finally, let $ k = 2$. The triple is then $(a,2a, 2a+\rho)$. 

We   know (more generally, independently of the value of $k$) that $\rho \equiv 1 \pmod a$. Modulo $y = 2a + \rho$, we have $a + 2a \leq 2a + \rho + 1$, so $\rho \geq a+1$. If $\rho > a+1$, then $\rho \geq 2a + 1$, and so modulo $y$ with $t =2$, we obtain $2a + 4a \leq 2a + \rho + 2$, so that $\rho \geq 4a-2$. Thus if $a > 2$, we have $\rho \geq 4a + 1$, so we can now take $t = 3$, and obtain $3a + 6a \leq 2a + \rho + 3$, that is, $\rho \geq 7a -3$. But then $y \geq 9a -3 > 6a = 3x$. However, $y \leq 3x$ following from the original inequality. We are thus reduce to the case that $a = 2 = k$. The triple is thus $(2,4,4+\rho)$. Of course, $\rho$ is odd. When $\rho = 1$, we obtain the triple $\boxed{(2,4,5)}$. If $\rho \neq 1$, then $\rho \geq 3$, so that the triple is in the form $(2,4,5 + 2s)$ where $s\geq 1$. Here $t_0 = 4$, so $ y \notin[ 8,\infty)$, leaving $7$ as the only possible value for $y$. But $(2,4,7)$ is not asymptotically hollow (modulo $7$ and with $t = 3$, we have $6 + 5 > 7+3$).\qed 
 
 \Pf (of (A) Lemma \Sfou) We can assume $a$ does not divide $x$; since $x \geq a^2$, we
thus have $a \geq x^2 +1$. With $a =r$, we have $\brm x,ty. \leq x - t(r-1)$ for
all $t \leq x/2$. Hence Theorem \Aone (in the Appendix)  applies, so that $y \equiv
1,x-a, x-a+1 \pmod x$. We go through the various possibilities. 

 \noindent $y \equiv 1 \pmod x$. Then $y = sx +1$ for some integer $s \geq 1$,
but since $a > 2$, each of these lies in the proscribed interval $[sx,
s(x+a-1))$. 

\noindent $y \equiv x-a, x-(a-1) \pmod x$. Then $y = sx-a$ or $y = sx-a+1$ with
$s \geq 2$. But $2/k + 1/(a-1) \leq 1$ if $a \geq 4$ or $k \geq 4$, in which
case $y <2x$ ($y = 2x$ is always a proscribed value), or $a = 3$ and $k = 3$,
and then $y <13 x/6$. In the latter case, we have proscribed interval
$[2x, 2x + 2a-3]$, $x = 9 + \delta \leq 11$, so $2a-3 = 3 > x/4 > x/6$. Hence in
 this case as well, we have $y < 2x$. Thus $ y = 2x - a$ or $2x -a +1$. 
 
 Now we consider the proscribed interval, $[(2-1/k)x, (2-1/k) (x+a-1))$
(obtained by setting $t = 2k-1$ in $[xt/k, (x+a-1)t/k)$). We claim that
$$
 \(2 -\frac 1k \)x \leq 2x-a, \qquad2x-a+1 < \(2 -\frac 1k \)(x + a -1),
$$
that is, either choice for $y$ is in a proscribed interval, hence the triple is
not \ah. The left inequality, $2x-a \geq (2-1/k)x$ is equivalent to $x/k \geq
a$, which of course is true (since $k = \flo{x/a}$). The right (strict)
inequality, $\(2 -1/k \)(x + a -1) > 2x-a+1$ is equivalent to $x/k <
(3-1/k)(a-1)$. Writing $x = ka + \delta \leq (k+1)a -1$, the left side is at
most  $a(1+ 1/k) - 1/k = (1+ 1/k)(a-1) + 1$, so the difference (right side minus
left side) is at least $2(a-1) -1 >0$. 
\qed 

  \Lem (B) Lemma \Ssix. No \ah\ triples exist with $a > k \geq 3 $. 
 
 Here we show that various proscribed intervals overlap, so their union contains
a large interval.   The following results are routine, but require a little care. 

\Lem Lemma \Ssev. Suppose that $1 \leq s \leq k-2$ and $a > k$.  Then 
$$
 \( 2- \frac s{k-1}\) (x+a-1) > \( 2- \frac sk\) x.
$$
 \Pf The inequality is equivalent to 
$$\eqalign{
\( 2- \frac s{k-1}\)(a-1) &> \frac{sx}{k(k-1)}; \quad{\text{as $x \leq
(k+1)a -1$, sufficient is }} \cr
 (2k^2 - 2(1+s)k - s)(a-1) & > sk \cr
}$$
The coefficient of $s$ in the left side minus that of the right side is $-(2k+1)(a-1) -
k < 0$; hence the difference is decreasing in $s$. It thus suffices to verify
the inequality when $s= k-2$. At $s = k-2$, the left side becomes $(k+2)(a-1)$,
and the right side is $k^2-2k$. Hence sufficient is 
$$
a > 1 + \frac{k^2 -2k}{k+2} = k-3 + \frac{8}{k+2}.
$$
Since $a > k$, sufficient is that $8/(k+2) \leq 3$, for which $k \geq 1$  is
enough.\qed

\Lem Lemma \Seig. Suppose that $ 2 \leq s  $ and  any one of the following pairs
of conditions hold: 
\item{(a)} $k \geq 6$ and $a \geq 3$; 
\item{(b)}$k=5$ and $a \geq 4$; 
\item{(c)} $k =4 $ and $a \geq 5$.
{\par} \noindent Then 
$$
 \( 2- \frac s{k}\) (x+a-1) > \( 2- \frac{s-1}{k-1}\) x.
$$
 \Pf This is equivalent to 
$$\eqalign{
\( 2- \frac s{k}\)(a-1) &> \( \frac sk - \frac{s-1}{k-1}\)x;\quad
\text{  $\times k(k-1)$; \text{as $x \leq
(k+1)a -1$, sufficient is},}\cr
 (k^2 - 2k - s)(a-1)  & > (k+\delta)(k-s).\cr
}$$
The right side is bounded above by $(k+ a-1)(k-s)$, so sufficient is 
$$
 (k^2 - 3k)(a-1) > k^2 -sk.
$$
 The right side is decreasing in $s$ (and the left side has no $s$ term), so it
is sufficient to verify the last inequality in the case that $s = 1$, that is, 
$$\eqalign{
 a- 1& > \frac{k^2 - k} {k^2 - 3k} = 1 + \frac{2k}{k^2 - 3k}; \quad \text{equivalently,}\cr 
 a &> 2 + \frac{2k}{k^2 - 3k}  = 2 + \frac 2{k-3}\cr
}$$ 

 If $k = 4$, $ a \geq 5$ is sufficient; if $k = 5$, then $a \geq 4$ suffices, and if $k \geq
6$, then $a \geq 3$ is sufficient. \qed

 \comment
 View the difference as a linear polynomial in $s$; the coefficient of $s$ is
$-(2k+1)(a-1) -k < 0$. Hence the difference is decreasing in $s$. It thus
suffices to verify the last inequality when $s = k-2$. Substituting, the left
side becomes  $(2k^2 - 2(k-1)k - (k-2))(a-1) = (k+2)(a-1)$. The right side is
$k^2 - 2k$. Thus sufficient is 
$$\eqalign{
 a-1 &> \frac{k^2 - 2k}(k+2) = k - 4 + \frac{8}{k+2}; \quad\text{equivalently,}
\cr 
 a > k-3 + \frac{8}{k+2}\cr
}$$
 Since $a > k$, sufficient for the last displayed inequality to hold is that
$8/(k+2) \leq 3$, for which $k \geq 1$ is sufficient. \qed
 \endcomment  
 
 Suppose that $k \geq 4$. Label some of the proscriptive intervals,  $I_i = [(2- i/k) x, (2- i/k)(x + a-1))$ for $1 \leq i \leq k
$, and   $J_i = [(2- i/(k-1)) x, (2- i/(k-1))(x + a-1) )$ for $1 \leq i \leq k-1
$. The  preceding two lemmas (when they apply) yield
$I_{i+1} \cap J_{i} \neq \emptyset$ and $J_s \cap I_{s} \neq\emptyset$. It
follows easily that the union, 
$$
I_{k-1} \cup J_{k-2} \cup I_{k-2} \cup \cdots\cup I_2 \cup J_1 \cup I_1, \quad\text{is the interval $\left[\(1 + \frac1{k }\)x,  \(2-\frac1k\)(x + a -1)\)$}.
$$
(To see this, read the union from the right to left, that is, starting with $I_1$.) Combined with Lemma \Sthi(b), and the fact that the very first proscribed in interval is $[x, x+a-1)$, we obtain the following.

 \Lem Lemma \Snin. If $a > k \geq 3$, then $ x+a-1 \leq  y < \(1 + \frac1{k}\)x$.

\Lem Lemma \Sten.  Suppose that $a \geq 3$. If   $k\geq 3$, or if $k = 2$ and $a\geq 4$, then $y\geq (1+ 1/k)(x+a-1)$.
 
 \Pf  We have the first proscriptive interval, $[x, x+a-1)$, so $ y \geq x+a-1$. From the proscriptive interval $[(1+ 1/k)x, (1+ 1/k)(x+a-1))$, it suffices to show that $\rho \not\in [x+a-1,(1+ 1/k)x]$. If $\rho$ does belong to this last interval, then $y = x + \rho$ with $a-1 \leq \rho \leq x/k$. Then $k(a-1) \leq k\rho \leq x$. Hence $\brm x,ky. = \brm x,k\rho. = k\rho$; since $\brm x,ka. = ka$, we have, with $t= k$ (which is obviously bounded above by $x/2$), applied to $(a,x,y)$, we have $ka + k\rho \leq k + x = k + ka + \delta$. Hence $k\rho \leq k+ \delta$, that is, $\rho \leq \delta/k + 1$. 

 As $\delta \leq a-1 \leq \rho$, we have $a-1 \leq (a-1) /k + 1$, that is, $k \leq (a-1)/(a-2) = 1 + 1/(a-2)$. The right side is at most $2$, so if $k \geq 3$, we are done; if $a \geq 4$, then the right side is strictly less than $2$, so we obtain a contradiction if $k = 2$. \qed


 \Lem Lemma \Sele. (a)(trivial) $y \not\in [2x-a+2, 2x-1]$.
\item{(b)} If $k \geq 2$ and $a \geq 3$, then $y \not\in [(2-1/k)x, 2x + 2a-3]$. 

 \Pf (a) Write $\rho = x-u$ where $1 \leq u$. With $t=1$ and modulo $x$, we have $a + x-u \leq 1 + x$, and thus $u \geq a-1$, so $y = 2x-u \leq 2x-a+1$. 

\noindent (b)  By (a) and the two proscriptive intervals, $[x(2-1/k), (x+a-1)(2-1/k)$ and $[2x, 2x + 2a-2)$, it suffices to show 
$$
 (2-1/k) x \leq 2x - a + 2 < (2-1/k)(x+ a-1).
$$
 The left inequality is equivalent to $x/k \geq a-2$, which is trivial since $x \geq ka$. The right inequality boils down (after simplifying and multiplying by $k$)  to $x < (3k-1)(a-1) -k$. Since $x \leq (k+1)a-1 = (k+1)(a-1) +k$, it suffices to show $(2k-2)(a-1) > 2k$. As $k > 1$, sufficient is thus $a -1 > k/(k-1)$, so $a \geq 3$ is enough. \qed

\Pf (of (B) Lemma  \Ssix)  Since $2/k + 1/(a-1) \leq 1$, so $ y\leq 2x$. By Lemmas \Snin--\Sele, there are no choices for $y$. \qed

 \Lem (C) Lemma \Stwe. The only asymptotically hollow triple with  $k = 2$ and $a \geq 3$ is $(3,8,10)$.

\Pf First, for
  $k = 2$ and $ a\geq 4$; we show no triples exist.  Candidate triples are of the form $(a,2a + \delta, 2a + \rho + \delta)$. Then $x/(a-1) = 2 +( \delta+2)/(a-1)$, and thus $\ceil{x/(a-1)} \in \brcs{3,4}$. Hence $ y \leq 4(2a+ \delta)/2 = 2x$, and since $2x$ belongs to a proscriptive interval, $y \leq 2x-1$. 

 By Lemmas \Sten\ and \Sele,  with $k =2$, we have $3(x+a-1)/2 \leq y < 3x/2$, a contradiction.

This leaves the possibility that  $k = 2$ and $a =3$; we must show that $(3,8,10)$ is the only one.   The triple is $(3, 7 + \epsilon, 7 + \epsilon + \rho)$ ($\epsilon \in \brcs{0,1}$). Here $\ceil{x/2} = 4$, so $y \leq  2x$, and thus, $y \leq 2x-1$. By Lemma \Sele(b), $y < 3x/2$, so $\rho < x/2$; since $a > 2$, $\rho \geq 2$ (from the first  proscribed interval $[x, x+1]$).  So $2 \leq \rho< x/2$; since $x$ is $7$ or $8$, we can only have $\rho\in \brcs{2,3}$. 
 
 With $\epsilon =0$, $x = 7$, and we are left with $y \in \brcs{9,10}$. The triples are thus $(3,7,9)$ and $(3,7,10)$. Modulo $7$ and with $t = 2$, we see neither is asymptotically hollow. 
 
 With $\epsilon = 1$, we have $x = 8$ and the corresponding triples are $\boxed{(3,8,10)}$ and  $(3,8,11)$. The first is asymptotically hollow, while the other is not  (modulo $8$, set $t = 2$).  \qed 

\Lem (D) Lemma \Sthr. The  only asymptotically hollow triples  with $k \geq 2$ and $a= 2$ are  $(2,x,x+1)$ ($x \geq 4$) and $(2,5,9)$. 

\Pf By Lemma \Sfiv, we can assume that $x$ is odd, so $x = 2k+1$. Here $x \geq a^2 = 4$, so by Theorem \Aone, $y \equiv \pm1, -2 \pmod x$.

First, consider the case that $k \geq 3$; we show no triples exist other than $(2,x,x+1)$.  Here $x \geq 7$.
We also have $y < \ceil{x/(a-1)}\cdot x/k = 4k + 4 + 1/k$ by Lemma \Sthi(c), so $y \leq 4k+4 = 2x + 2$. Hence $y$ can only be one of $x+1, 2x-2, 2x-1$, or $2x+1$. The first yields $\boxed{(2,x,x+1)}$. We eliminate $2x+1 \in [2x,2x + 2)$, a proscriptive interval. The proscriptive interval $[(2-1/k)x, (2-1/k)(x+1))$ contains both $2x-2$ and $2x-1$: $(2-1/k)x = 2x- 2 - 1/k \leq 2x-2$, and $(2- 1/k)(x+1) = 2x+2 -2-2/k > 2x-1 $.

\noindent$k = a = 2$. We must have $x = 5$, so the triples are of the form $(2,5, 5+\rho)$. The proscribed intervals are $[5,6)$, $[15/2, 9)$, $[10,12)$, $[25/2, 15)$, and the union of the remaining ones includes $[15, \infty)$.  This leaves possible values for $y $, $6,7,9,12$. With $ y =6$, we obtain $\boxed{(2,5,6)}$. We also have that $ y\equiv \pm1,-2 \pmod 5$, eliminating $7$ and $12$. The only remaining one is $y = 9$, yielding $\boxed{(2,5,9)}$.   \qed   

 \comment
 Case $k \geq 3$, $a \geq 4$. No triples.  Here $y \leq (2/k + 1/(a-1)) x \leq 2x$, that is, $y \leq 2x-1$. By Lemma B (b), $y < (2-1/k)x$; by Lemma xxx, $y \leq (1+ 1/k)x$, and this contradicts Lemma A. \qed  

\endcomment
 
 \Lem (E) Lemma \Sftn. With $k =1$, the only asymptotically hollow triples are $(2,2,3)$, $(2,3,4)$, $(2,3,5)$, $(2,3,8)$, $(3,4,6)$, $(3,5,7)$, $(3,5,8)$, $(4,6,9)$, $(4,7,10)$, $(5,8,12)$ and $(6,10,15)$. 

\Pf Here $x = a + \delta$. 

\item{(i)}First, we consider the  case wherein  $ a\geq 3$, and $\delta = a-1$. The triple is then $(a, 2a-1, 2a-1 + \rho)$. 

We have $y \leq (2 + 1/(a-1)) x <5 x/2$, so $\rho < 3x/2$.
 If $\rho \leq x$, then modulo $x$, we obtain $a + \rho \leq x+1 = 2a$, so $\rho \leq  a$. We can rule out $\rho  = 1$ (since $x+1$ is in a proscriptive interval), so $\rho > 1$ and thus $\rho-1 > 0$, and in particular, $\brm a,2a-1+ \rho. = \brm a,\rho-1. = \rho-1$. We thus obtain $a-1 + \rho-1 \leq 1 + a$, and so $\rho \leq 3$. This leads to only two possibilities $\rho = 2$ or $\rho = 3$, and the resulting triples are $(a, 2a-1, 2a+1)$ and $(a,2a-1, 2a+2)$. Modulo $y$ with $t = 1$, we obtain respectively $a + 2a-1 \leq 1 + 2a+1$ and $a + 2a-1\leq 2a+3$, that is, $a \leq 3$ and $a \leq 4$ respectively. 
 
 We are thus reduced to the following triples. 

 \noindent $\rho = 2$:  $\boxed{(3,5,7)}$, $(4,7, 9)$.

 \noindent $\rho = 3$:  $\boxed{(3,5,8)}$, $\boxed{(4, 7, 10)}$. 
 
 We eliminate  $(4,7,9)$ (modulo $9$ and $t=1$).  

 
 \item{(ii)} Now suppose that  $a \geq 3$, and  $\delta \leq a-2$.  

Then $(a + \delta)/(a-1) \leq 2$, so $ y < 2x$. Thus $x+ a-1 \leq y \leq 2x-a+1$, and the candidate triples are of the form $(a, a+ \delta, a+ \delta + \rho)$, with $a-1 \leq \rho \leq x-a+1 = \delta + 1$. As $\delta \leq a-2$, we must have $\delta = a-2$ and $\rho = a-1$. Hence the triples are of the form $(a, 2a-2, 3a-3)$. When $a = 3$, we obtain the known asymptotically hollow triple $\boxed{(3,4, 6)}$. If  $a \geq 4$, then modulo $a$ and with $t=1$, we have $a-2 + a-3 \leq 1+ a $, and thus $a\leq 6$. The resulting three triples (arising from $a = 4,5,6$) are all asymptotically hollow, 
$$
 \boxed{(4,6,9)},  \boxed{(5,8,12)},  \boxed{(6,10,15)}. 
 $$

 \item{(iii)}Finally, suppose  $k =1$ and $a = 2$. We have $x \in \brcs{2,3}$. 

With $(2,2,2+\rho)$, $\rho = 1$ yields $\boxed{(2,2,3)}$. $\rho$ must be odd (as  asymptotically hollow triples  have content one), and we see that $(2,2,2u+1)$ (with $u > 1$) is not asymptotically hollow: modulo $y = 2u+1$ and $t = u$, we have $4u \leq 1 + (2u+1)$, forcing $ u \leq 1$, a contradiction.  

 For $x = 3$ and corresponding triple $(2,3, 3+\rho)$, by Lemma \Sthi(c), $y <9$. Also each of $x, 2x, 2x+1 $ is in the union of the proscribed intervals $\Z \cap \([x,x+1) \cup [2x, 2x+ 2a-2) \)= \brcs{3,6,7}$, leaving possible values for $y$, $\brcs{4,5, 7,8} $. These yield  the only candidates 
 $$ 
 \boxed{(2,3,4)}, \boxed{(2,3,5)}, (2,3,7), \boxed{(2,3,8)}.
 $$
 The triple $(2,3,7)$ is eliminated by taking $t = 2$, modulo $7$. \qed

\SecT Appendix

Here we prove a result used in the proof of  the $n=4$ classification theorem (\Sone). Let $x,r$ be  integers exceeding $1$ with $x \geq 2r$. Define
$$\eqalign{ 
S(x,r) &:= \Set{z \in \brcs{1,2,\dots,x}}{\brm x,zt. \leq x - (r-1)t \text{ for all $t \leq x/r$}}\cr
S_0(x,r) &:= \Set{z \in \brcs{1,2,\dots,x}}{\brm x,zt. \leq x - (r-1)t \text{ for all $t \leq x/r$ \st $\brm x,tz. \neq x$}}\cr
}$$
We observe that $S(x,r) \subseteq S_0(x,r)$, and $z \not\in S_0(x,r)$  if and only if there exists an integer $t \leq x/r$  \st $x > \brm x,tz. > x - (r-1)t$.

\Lem Theorem \Aone. Let $x$ and $r$ be  integers exceeding $1$, with $x \geq r^2$.   
Then
$$S(x,r) = \cases \brcs{1,x-r,x-(r-1)} & \text{if $r$ does not divide $x$}\\ \brcs{1,x-(r-1)} & \text{if $r|x$ and, if $r =2$,  $x \equiv 2 \pmod 4$}\\ \brcs{1, \frac{x}2 -  1, x-1}& \text{if $r = 2$ and $x \equiv 0 \pmod 4$.}\\ \endcases $$

\Rmk The necessary and sufficient  condition for this result to be true  is likely $x \geq r\ceil{r/2}$  (and this is clearly necessary) rather than $x \geq r^2$.

The proof of   Theorem \Aone\ is rather tedious. In all cases, we proceed by eliminating candidate $z$s from $S = S(x,r)$. There are two general techniques, based on one idea.

\comment
The proof will be subdivided into cases. 

\item{(a)} If $z \in S(x,r)$, then $\gcd{(z,x)} <r$;
\item{(b)} $r|x$;
\item{(c)} $S(x,r) \cap [2,r] = \emptyset$;
\item{(d)} $S(x,r) \cap [r,3r/2] = \emptyset$;
\item{(e)} $S(x,r) \cap [r+4,x-2r] = \emptyset$;
\item{(f)} $S(x,r) \cap [x- (r-1)\sqrt{r+1},x-r-1] = \emptyset$;
\item{(g)} $S(x,r) \cap [x-r+2,x] = \emptyset$;
\item{(h)} any residual cases.

In each of (c--h), $r$ does not divide $x$. Several of these are quite easy. For the rest, there is a general technique used in many of the cases. 
Begin with $z$ not in the indicated set of three or two numbers. We wish to find a positive integer $s$ \st at least one of $\brm x,zt. > x - t $ for $t = \ceil{sx/z}$ or $\flo{sx/z}$, also requiring the corresponding value of $t$ to be at most $x/r$. This will yield a contradiction, but a combination of the two techniques will not yield the complete result. A few additional cases have to be considered (especially when $r = 2$). We analyze the possibilities. We first try $s = \flo{z/r}$, which we write as $(z - \varepsilon)/r$ (that is, $\varepsilon$ is the (usual) remainder modulo $r$, and write $x = vr + \delta$, again with $\delta$ the usual remainder modulo $r$. Then 
$$\eqalign{
\frac{sx}z & =\frac{ (z-\varepsilon) x}{rz}\cr  & = \frac x{r} \cdot \(1 - \frac{\epsilon} y \); \quad \text{in particular, $\flo{sx/y} \leq x/r$.}\cr 
 sx - (z-1) & \leq z\cdot \flo{sx/r} \leq sx;\quad \text{and thus}\cr 
\brm x,z\cdot \flo{sx/r}. & \geq x - (z-1). \quad \text{so if can show that }\cr z-1 &<  (r-1) \flo{sx/r}, \quad\text{then $z \not\in S$.} 
}$$ 
This  is effective when $z$ is smaller than slightly less than $x-(r-1)$ (and not one of the forbidden values). When $z > x-(r-1)$, we $y = x-z$, and now try to find $t$ of the form $\ceil{sx/y}$; in this case, the condition that $t \leq x/r$ becomes an issue, and the situation is more complicated (and involves the remainders $\epsilon$ and $\delta$; if $s = \flo{y/r}$ does not work (it could result in the corresponding candidate for $t$ being larger than $x/r$, we try $s= \flo{y/r} -1$. The combination of these two techniques leaves a few missing values of $z$ when $r$ or $x$ is small, and these have to be dealt with by ad hoc methods.   

\endcomment

\Lem Lemma  \Atwo. (a) If $\gcd \brcs{z,x} \geq r$, then $z \notin S(x,r)$. 
\item{(b)}  If $r|x$, then the conclusion of Theorem \Aone\ holds. 

\Pf (a) Suppose $a \geq r$ divides both $z$ and $x$.  Then $x/a \leq x/r$ and $z\cdot x/a \equiv 0 \pmod x$, and thus $\brm x,(z\cdot x/a). = x$, so $ z \not\in S(x,r)$. 

  \noindent (b) (i) $r|z$. Follows from (a). 

\noindent (ii) {\it $r \notdivide z$ and $z \not\equiv \pm 1 \pmod x$.} Say $z \equiv \delta \pmod r$ where $2 \leq \delta \leq r-1$. Set $t = x/r$.  Write $z = ur + \delta$. Then $z x/r = (ur+\delta)x/r = ux + \delta x/r \equiv \delta x/r \pmod x$. As $\delta > 1$, and  with $t = x/r$, we   have $\brm x,zt. = \delta x/r > x/r = x- (r-1)t$, so that $z \notin S(x,r)$.

\noindent (iii) {\it $r \notdivide z$ and $z \equiv \pm 1 \pmod x$.}

\noindent So assume $\delta = 1$, and for now, $r > 2$. We try to find a positive integer $w < x/r$, so that $t: = x/r - w$ will eliminate $z$ from $S$. 

\noindent {\it $x/r < z < x-(r-1)$}. There is a smallest positive integer $w$ \st $(w-1)z \leq x/r$ and $wz > x/r$ (if $z > x/r$, $w =1$). Then $z (x/r -w) \equiv x/r - wz \pmod x$, and thus $\brm x,yt. = x + x/r - wz$. This will exceed $x - (r-1)t = x-(r-1)(x/r -w)$ if and only if $wz - x/r < (r-1) (x/r -w)$; the latter is equivalent to  
$$x > w (z+r-1) \tag 1$$
  If $z > x/r$, then $w = 1$, so that (1) equivalent to $x > z+r-1$; thus if $x/r < z < x-(r-1)$, then $z \not\in S$. 

\noindent {\it $1 < z < x/r$.}  We have $r <  z+r-1 < x/r + r-1 < 2x/r$; since  $z+ r-1 \equiv 0 \pmod r$, we have that $z+ r-1 = (u+1)r$  with $u \geq 1$. If $w(y+r-1) \geq x$, then $x/r \leq w(u+1)< w(z/r) + 2 \leq z/r + x/r^2 + 2 < 2x/r^2 + 2$, that is, $x < 2x/r + 2r$. This is impossible, unless $r =2$, which requires a special case. 

\noindent{\it Remaining cases with $r > 2$.} If $z > x-(r-1)$ and $z \equiv 1 \pmod r$, then as $x \equiv 0 \pmod r$, $z$ cannot lie in $\brcs{x-r+2, \dots, x}$, and so this case is done. If $z = x/r$, then set $t = r$; as $r \leq x/r$ and $\brm x,x. = x$, we thus have $x/r \notin S$. Now the cases $r = 2$, $z$  is odd and less than $x/2$. 

\noindent {\it The case that $r=2$, $z$ is odd, and $1 < z < x/2$.} We apply the method outlined in the introduction to this appendix; the candidate for $t$ is $\flo{((z-1)/2)x/rz} = \flo {x(1- 1/z)/2}$, and we wish to verify that $z-1 < (r-1)t $. Since both sides are integers, this is equivalent to $z \leq (r-1)t$. As $t > x(1-1/z)/2 -1$, it is sufficient to show $z \leq (r-1)(x(1-1/z)/2 -1)$. This translates to a quadratic inequality in $z$,  
$$ z\(z- \(\frac x2 -1\)  \) + \frac x2 \leq 0. $$ 
The set of real values of $z$ for which this holds is a closed interval. At $z = 2$, the value of the quadratic is $-2 (\frac x2 -3) + \frac x2 = -\frac x2  + 6$; this is less than or equal zero when $x \geq 12$. At $z = x/2-3$, we obtain $(x/2-3) 2 + x/2 $, which is not positive when $ x \geq 12$. Thus $[2, x/2-3] \cap  S(x,2) = \emptyset$. At this point, we notice that if $4$ divides $x$, then $x/2 -2 $ is even, so the largest value we have to consider for $z$  is $x/2 -3$, and this completes the proof of the result in the case that $r=2 $ and $4|x$. 

Now suppose that $x \equiv 2 \pmod 4$, so we want to conclude that $z = x/2 - 2  \in S(x,2)$ (this will complete the proof for this case, since $x/2 -1$ is even). Set $t = x/2 - 2$, so that  
$$\eqalign{ 
zt &= \( \frac {x -4}2\)^2 = \frac{x^2 - 8x + 16}4 \cr & = x\frac{x-6}4 + 4 - \frac x2\cr & \equiv \frac x2 + 4 \pmod x \qquad\text{ since $x \equiv 2 \pmod 4$; so }\cr \brm x,zt. & = \frac x2 + 4 \quad \text{if $x \geq 8$}.  
}$$ 
On the other hand, $x - t = x - (x/2 - 2) = x/2 + 2 < x/2 +4$, so $x/2- 2 \notin S(x,2)$. \qed

From now on, we deal with the case that $r$ does not divide $x$. We have an easy   lemma of occasional use. 

\Lem Lemma \Athr.  With $t = \flo{x/z}$, sufficient for $\brm x,tz. > x - (r-1)t$ is that $\flo{x/z}(r-1)$ exceed the remainder of $x$ modulo $z$.

\Pf 
Write $x = kz + \eta$ where $0 \leq \eta< z$; we can exclude the case that $\eta = 0$ by Lemma \Atwo(a); and $x > z$, so $k \geq 1$. Then $x - (r-1)t = (kz + \eta) - (r-1)k = k(z -r+1) + \eta$. On the other hand, $tz =kz < x$, so $\brm x,tz. = tz$. Then the inequality holds if (and only if) $k z > k(z-r+1) + \eta$, that is, $k(r-1) > \eta$.

\Lem Lemma \Afou. Suppose that   $r$ does not divide $x$. 
\item{(a)}   $S_0(x,r) \cap [2,r] = \emptyset$. 
\item{(b)}   $S_0(x,r) \cap [2,3r/2] = \emptyset$. 

\Pf (a) Write $x = vr + \delta$ where $1 \leq \delta < r$. Set $t = v = \flo{x/r}$. Then $tz \leq tr < x < (r+1)t$. In particular, $\brm x,tz. = tz$. On the other hand, $x - (r-1)t = v+ \delta$. As $x > r^2$, we have $v \geq r$. Finally, $tz \geq 2v \geq v+r > v+\delta = x-(r-1)t$. 

\noindent In view of (a), if $r =2$, it is sufficient to show that  $3 \notin S(x,2)$; then $t = \flo{x/3}$ will work, as is easy to see from Lemma \Athr\ above. If $r =3$, we need only show $4 \notin S(x,3)$, and here again, $t = \flo{x/4}$ will do. If $r = 4$, only $5,6 \notin S(x,4)$ is required; but $x \geq 17$ entails $\flo{x/6} \geq 2$, and \wrt to Lemma above, $\eta < 6$ and $r = 4$. 

So we can assume $r \geq 5$. On the interval $r \leq z \leq 3r/2$, the corresponding $\varepsilon$ varies from $0$ to at most $\flo{r/2}$, and the ratio $\varepsilon/z$ is bounded above by $1/3$. This yields (in (2)), $t \geq \flo{2x/3r} > 2x/3r -1$. As $x/r > r$, we obtain $2x/3r -1 > 2r/3 -1$. Hence $(r-1)t > (2r/3 -1)(r-1)$. As $z \leq 3r/2$, it suffices to verify $3r/2 \leq (2r/3-1)(r-1) $. The quadratic $2r^2/3 - r(3/2 + 2/3 +1) +1 $ is positive for $r \geq 5$. 
\qed

\comment
Now we have a general construction method, variants of which will be used in subsequent cases. Write $x = vr + \delta$ where $1 \leq \delta < r$. Pick $z \in \brcs{2,\dots,x-1}$, and set $s = \flo{z/r}$, so that $z = sr + \varepsilon$ where $0 \leq \varepsilon < r-1$. The candidate for $t$ to eliminate $z$ (that is, to show it does not belong to $S(x,r)$ is $t = \flo{sx/z}$ (this does not always work; but minor variations of this can apply). From $sx/z = (x/r)(1- \varepsilon/z)$, we see that this choice of $t$ satisfies $x \leq x/r$. 

Now $\flo{sx/z}z = sx - \alpha$, where $0 \leq \alpha < z$. Hence $\brm x,tz. \geq x-\alpha \geq x - (z-1)$. Then $\brm x,tz. > x-(r-1)t$ will occur if $z-1 < (r-1)\flo{sx/z}$. As both sides are integers, $z \leq (r-1)t$ is sufficient.

We have 
$$\eqalign{ \frac {sx}z & = \frac{(vr+ \delta)\frac{z-\varepsilon}r}{z} \cr
& = v\(1 - \frac{\varepsilon}z\) + \frac {\delta}{r}\(1 - \frac{\varepsilon}z\)\cr 
& = \frac xr \(1 - \frac{\varepsilon}z\). \cr
}\tag 2$$
At this stage, we have various options, exploiting $x/r > r$. For example, if $\varepsilon \leq z/2$, then $sx/r \geq x/2r $, and thus $\flo{sx/z} \geq \flo{r/2} \geq (r-1)/2$. Or, we can just use $\flo{sx/z} > (x/r)(1-\varepsilon/z) -1$. 

\endcomment
 
Now we have a  reflection principle, but only for $S_0(x,r)$.  

 \Lem Lemma \Afiv. Suppose that $z \leq x-r$ and $z \not\in S(x,r)$. 
Then $x - (r-1) -z \not\in S_0(x,r)$.   In addition, if $(x,z) = 1$, then $z\not\in S(x,r) $ entails $z \not\in S_0(x,r)$.

  \Pf There exist nonnegative integers $t$ and $b$ with $t \leq x/r$ \st 
$x > tz - bx > x - t(r-1)$.

 Then we have 
$$\eqalign{
 t(x-z -(r-1)) & = tx -tz - t(r-1) \cr
 & > tx - (b+1)x -t(r-1) \qquad \text{so}\cr
 t(x-z-(r-1)) - (t-b-2) x & > x - t(r-1); \qquad {\text{moreover,}} \cr
t(x-z-(r-1)) &= tx - tz - t(r-1) \cr 
 &< tx - bx -x + t(r-1) - t(r-1), \qquad\text{and thus,}\cr 
 t(x-z - (r-1)) - (t-b-2) & < x .\cr
}$$

If $x = tz-bx$, that is, $tz = (b+1)x$, since $(z,x= 1)$, we have $z$ divides $b+1$. But $t \leq x/r$ entails $z \geq (b+1)x/r > (b+1)$, a contradiction. \qed 
 
 The following  easy special cases will reduce the number of cases to be dealt with later. 
 
 \Lem Lemma \Asix. (a) If $x \geq 7$, then $[(x+3)/2-r, (x-(r-1))/2]\cap S_0(x,r) = \emptyset$. 
\item{(d)} If $x-r+2 \leq z \leq   x$, then $z \not\in S(x,r)$.
\item{(e)} ($r=2$). Let $x$ be an odd integer exceeding 20, and let $p \in \brcs{3,5,7}$. Then $(x-p)/2 \not\in S(x,2)$. {\par}
\item{(f)} Suppose $x \geq r^2$ and $1 \leq z \leq r$ is an integer of the form $z = kr+l$ where $k$ and $l$ are nonnegative integers \st $z = kr +l$. If $k+2 \leq l \leq r$, then $z \not\in S(x,r)$.  \item{(g)} Suppose that $x \geq 10$. Let $z = r+ \epsilon$ where $\epsilon \in \brcs{1,2}$. Then $z \not\in  S(x,r)$. 

 \Pf (a) Set $t =2$ and let $z = (x-a)/2$. Then $\brm x,tz. = x-a$ where $a \in [r-1,2r-3]$, whereas $x - (r-1)t = x-2(r-1)$. Then $a < 2(r-1)$ is sufficient for $z \not\in S_0(x,r)$.
 
\comment
 (b) Set $z = \flo {\frac { x-(r+1)}2}$. First, set $t=2$.  If $x-r$ is odd, $tz = x-r-1 > x- 2(r-1)$ if $r \geq 4$; and if $x-r$ is even, then $tz = x-r-2$, and this exceeds $x-2(r-1)$ if $r \geq 5$. 
 Now suppose that $r = 3$ and $x-3 = x-r$ is odd, so that $z = (x-4)/2$.

 \noindent (b) If $x$ is even, we are done by xxx. Otherwise, we have $z = (x-3)/2$ is an integer. We wish to find a positive integer $s$ \st with $t = 2s+1$, we have $t \leq x/2$ and  $\brm x,tz. > t$. The first constraint is that $s \leq (x-3)/4$. 

 To determine what the  second is equivalent to, we note 
 $$
   (2s+1)\frac{x-3}2 = s(x-3) + \frac {x-3}2 \equiv  \frac{x-1}2 - 3s \pmod x,
$$
  so that if $s \geq (x-1)/6$, then $\brm x,tz. = (x-1)/2 - 3s + x$ (assuming $s < (x-1)/2$, which would be a consequence of $s \leq (x-4)/3$). So we want $(x-1)/2 - 3s + x > x - 2s-1$, which boils down to $s < (x+1)/2$, a consequence of the first constraint. 
 
 Thus, if  $(x-1)/6 \leq s \leq  (x-3)/4$, the corresponding $t$ will yield the result.  It therefore suffices to find an integer in the interval $[(x-1)/6,  (x-3)/4]$. It  has length $(x-7)/12 $; this is at least one (and thus the interval contains an integer) if $x \geq 19$. If either $x \equiv 3 \pmod 4$ or $x \equiv 1 \pmod 6$, then one of the endpoints is an integer. The only remaining odd value of $x$ between $13$ and $19$ is $x= 17$, for which the interval does contain an integer. 

 \noindent (c) Set $t = \flo {x/r}$. 

 \endcomment 

 \noindent (d) Set $t =1$.

\noindent (e) Set $k = \ceil{x/p}$, and let $s = kp - x$; obviously, $0 \leq s \leq p-1$. Then $k(x-p)/2 = (k-1) x/ 2 - s/2 $. 

 If $k$ is odd,  set $t = k$. Then $\brm x,t(x-p)/2. = x - s/2$. Then $s/2 < t =( x+s)/p$ provided  $s (p/2 - 1) < x$. This is implied by $(p-1) (p/2 -1)< x$. The biggest value attainable for the left side is $15$, so $x \geq 17$ is sufficient. Obviously $t \leq x/p + 1$, and this is bounded above by $x/2$ (as $ \geq 3$ and $x > 19$).
 
 Now suppose that $k $ is even. Set $t = k+1 $, so that $tp = x + p + s$. Then $t(x-p)/2 = tx/2 - tp/2  = (t-1)x/2- s/2 - p/2 $. Since $k = t-1$ is even, $\brm x,t(x-p)/2 . = x-(p+s)/2$. Also $t \leq x/2$ is implied by $x/p + 2 - 1/p < x/2$, that is, $(p-2)x > 4p-2 $, and this is true for $x \geq 11$. So it remains to determine when $(p+s )/2 < t $. This boils down to $x > p(p+s-2)/2 - s$. 
 
 Worst case occurs when $p=7 $ and $s =6$, so that $x > 6(13)/2 -6 = 33$ is sufficient. This leaves the values $x =33, 31, 29, \dots , 21$. When $p= 7$, $k = \ceil{x/p}$, and this is odd for $x = 33, 31, 29$.  For $x = 27$, we have $k= 4$, $s = 1$, and $t= 5$; $p(p+s-2)/2 - s = 7\cdot 6/21 - 2 = 19 < 27$. If $x = 25$, $s = 3$ and $p(p + s -2)/2 - s = 7\cdot8/2 - 2 = 24 < 25$.  When $x = 23$, $s = 5$. Here $(x-p)/2 = 8$; here we take $t = 8 \leq 11$ instead. When $x = 21$,  here $z = 7$, and we just take $t = 3$. 
 
 When $ p =5$, $x > 5\cdot 7/2 - 4 < 14$ and there is no problem (nor for $ p =3$).

\noindent (f) Set $t = \flo{x/r}$, that is, $t = (x-\delta)/r$ where $0 \leq \delta \leq r-1$.  Then $tz = kx - k\delta + l(x-\delta)/r$. Modulo $x$, this is equivalent to $u:= l(x-\delta)/r - k\delta$. Obviously $u \leq x$; we claim that $u > x - (r-1)t$  (which also entails, in this case, that $u = \brm x,tz.$.

 The inequality $u > x-(r-1)t$ is equivalent to 
$$
 x  > \delta\( 1 + \frac{(k+1)r}{l-1}\). 
$$
 But the right side is bounded above by $(r-1)(1 + r) = r^2-1$, and the result follows.

\noindent (g) Set $t = \flo {x/z} = (x-\phi)/z$. Then $t \leq x/r$ and $tz = x-\phi$, and thus $\brm x,tz. = x- \phi$. So it suffices to show that $\phi < (r-1)t$. This is equivalent to 
$$
\phi < \frac{r-1}{r+ \epsilon} (x-\phi);
$$
this is in turn equivalent to $(r-1)x > (\epsilon + 1)\phi$. Since $\phi \leq  r + \epsilon -1$, sufficient is $x > 3(r+\epsilon -1)/(r-1)$. The right  side is at most $9$, and we are done. 
\qed 

\comment 
  \noindent (e) Write $x = kz + \eta$, with $0 \leq \eta \leq z-1$. If $\eta < (z-1)/2$, set $t = k = \flo{x/z} = (x-\eta)/z$. As $z \geq r$, it follows that $t \leq x/r$. We have $tz = x-eta$, and thus $\brm x,tz.= x-\eta$. On the other hand, $x - t(r-1) = x- (x-\eta)(r-1)/z$. Sufficient is thus $\eta < (x-\eta)(r-1)/z$, which is equivalent to $\(1+(r-1)/z \)\eta   < x$. Since $\eta \leq (z-1)/2$ and $r-1 < z$, the left side is less than $z-1$, and we are done.

 Now assume that $\eta \geq (z-1)/2$.  We have $2x \equiv 2\eta \pmod z$, and we try $t = \flo{2x/z}$. Now $2\eta \geq z-1 \geq r$, and  
 
\endcomment

 \Lem Corollary \Asev. Suppose $x \geq r^2 \vee 10$. Then $\brcs{2,3, \dots, 2r} \cap S_0(x,r) = \emptyset$.  
 
 \Pf For $2 \leq z \leq r$, the result follows from Lemma \Asix(f) with $k = 0$; for $r+3 \leq z \leq 2r$, it follows when $k = 1$. Finally, for $z = r+1$ or $r+2$, it follows from Lemma \Asix(g). \qed


\comment
\Lem Lemma \Aeig. Suppose that  $r$ does not divide $x$. Let $F(z) = (z-1)^2 + r(z-1) - (r-1)x $. If $F(z) < 0$, then $z \not\in S(x,r)$. 

\Pf Write $x = lz + \eta$, where $\eta \leq z-1$. By Lemma \Athr,   sufficient   that $z \notin S(x,r)$ is that $\eta <\flo{x/z}\cdot(r-1) $; this amounts to $(x-\eta)(r-1)/z > \eta$, which we rewrite as $\eta (z + r-1) < x(r-1)$. Sufficient (using $\eta \leq z-1$) is thus 
$$\eqalign{
(z-1) (z-1+r) -(r-1)x &< 0; \qquad \text{equivalently,}\cr
(z-1)^2 + r(z-1) - (r-1)x & < 0.
}$$\qed

\Lem Corollary \Anin.   If $r$ does not divide $x$, then $S(x,r) \cap [0,\sqrt{(r-1)x} -r/2 + 1) = \emptyset$.

\Pf Treating $z-1$ as a real variable, the set where the quadratic $Y^2 + rY - (r-1)x$ is negative constitutes   an open interval containing $0$. The positive zero, $z_0-1$, is given as 
$$\eqalign{
z_0-1 & = \frac 12\(-r + \sqrt{r^2 + 4(r-1)x} \)\cr
& = -\frac r2 + \sqrt{(r-1)x} \( 1 + \frac {r^2}{4(r-1)x}\)^{1/2}; \quad \text{using $r^2 < 4(r-1)x$, i.e., $4x > r+1 + 1/(r-1)$,}\cr
& = -\frac r2 + \sqrt{(r-1)x} + \frac{r^2}{8\sqrt{(r-1)x}} - \frac{r^4}{128((r-1)x)^{3/2}} + \dots.\cr
}$$
The series is alternating and the absolute value of successive terms is decreasing. Hence $z_0 \geq \sqrt{(r-1)x} -r/2 + 1$, and $1 < z < z_0$ is sufficient. \qed

\endcomment

\Lem Lemma  \Aten. Suppose $x \geq r^2 +1$. Set $f(z) = z^2 - (r-1)\( \frac xr-1\)z + (r-1)x$. If $f(z) \leq 0$, then $z \not\in S_0(x,r)$.

\Pf Given the integer $z$ in $[r+4, x-2r]$, we wish to find integers $a = t$ and $b$ \st the following hold: 
\item{(i)} $x > az-bx > x-(r-1)a$; 
\item{(ii)} $1 \leq a \leq x/r$. 

Equation (i) can be rewritten as 
$$
\frac{(b+1)x}z > a  >\frac{(b+1)x}{z+r-1} \tag 3
$$
The  open interval $I:=\((b+1)x/(z+r-1),(b+1)x/z\right)$  will contain an integer if its length exceeds one; so we will require the following inequality to be satisfied by $b$:
$$
(b+1)x(r-1) > z(z+r-1). 
$$
If this holds, we will choose any integer $a \in I$. To guarantee that the resulting choice satisfies (ii), we also impose the condition that $(b+1)/z \leq 1/r$ (from (3), this will imply $a \leq x/r$). So we are looking for the positive integer $b+1$ to be in the interval
$$
\frac{z(z+r-1)}{x(r-1)} < b+1 \leq \frac{z}{r}  \tag {$3'$}
$$
To find such $b$, it is sufficient that $z \geq r$ (we only require $b \geq 0$) and the interval has length at least $1$. The length of the interval is $z/r - (z/x)(z+r-1)/(r-1)$, and for this to be at least $1$ is equivalent to the inequality, 

$$\eqalign{
f(z):= z^2 - (r-1)\( \frac xr-1\)z + (r-1)x  &\leq 0\cr
}$$\qed

Here $z$ is allowed to vary over $(r,x]$; viewing $z$ as a real variable, the set of 
$z$ at which this  quadratic is non-positive is a closed interval.

\comment
 
\Lem Lemma \Aele. Suppose $x \geq (r^2 + 1) \vee 13$. 
\item{(a)} Then $F(2r+1) < 0$
\item{(b)} If $x \leq 10(r-1)$, then $F((x/2 +1 - r) < 0$. 

\Pf  (a) With $F(z) = (z-1)^2 +  (z-1)r - 2(r-1)x$ 
 we have $F(2r + 1) = 4r^2 + 2r(r-1) -2(r-1)x$, so is negative iff $x > (3r^2 - r)/(r-1) = 3r + 2 + 2/(r-1)$. When $x \geq r^2 + 1$, this is automatic if $r \geq 4$. For $r =2, 3$ respectively, the corresponding conditions are $x \geq 11, 13$ respectively (and for these, just a few cases,e.g., $x = 11$ and $r =2$, just $z = 5$ has to be checked, since $\brcs{2,3,4}$ are already done. 

  \noindent (b)  
 $$\eqalign{
 4F\( \frac{x}{2} + 1 -r\)& =  \(x- 2r\)^2 +2\( x  -2r\)r - 8(r-1) x \cr 
 & = x^2 - x (10r-9) .\cr
 }$$
 In particular, $F( \frac{x}{2} +1 -r) < 0$ if $x \leq 10r -9$.
\qed
 
 In view of the preceding results, if $x \geq r^2 +1$, either $r \geq 4$ or $x > 3r + 2 + 2/(r-1)$, and $x \leq 10(r-1)$, then the result holds. The middle condition can be dropped if we verify the few remaining cases with $r = 2,3$. This means that we can now assume that $x \geq 10r-9$ (and for $r = 2,3$, somewhat larger numbers). 
 Now set $f (z) = z^2 - (r-1)\( \frac xr-1\)z + (r-1)x $. Then $rf(2r+1) =  r(4r^2 + 4r + 1) - (r-1)(x-r)(2r+1) + r(r-1)x$, which simplifies to 
$$\eqalign{
 x(r^2-r -(2r^2-2r + r -2) + &4r^3 + 4r^2 + r +(2r+1)r (r-1) \cr
  = x(-r^2 + 2) &+ 6r^2 + 3r^2. \cr 
}$$
 This will be negative if $$
  x > \frac{6r^3 + 3r^2}{r^2-2} = 6 r + 3 + \frac{12r + 6}{r^2-2}. 
$$
 Compare this with $10r - 9$; the difference is $3r -6 - (12r + 6)/(r^2 - 2)$. At $r = 4$, the difference is $6- 54/14 > 0$, and similarly, for $r \geq 5$, $10r-9 $ is larger than the displayed value, and thus the $f$ value is negative.  (can we do better: $f(\sqrt{2(r-1)x}- r/2)$?)

\endcomment

\Lem Lemma \Atwe. Suppose that $x \geq (r^2 + 1)\vee (10r-8)$, and $r > 3$.  If $x \geq 5r$, then $f(2r+1) \leq 0$; if $x \geq 10r-8$ and $r \geq 4$, then $f(x/2 +1 -r) \leq 0$.

\Pf Recall that $rf (z) =  rz^2 - (r-1)(x-r)z + r(r-1)x$.  
Since the set of $z \in \R$ \st $f(z)< 0$ is a (possibly  empty) open interval, it suffices to show that  $rf((x+2)/2 - r) < 0$ and $rf(2r+1) < 0$.

To show $rf((x+2)/2 - r) < 0$, we observe that   $\frac{\partial(rf)}{\partial x} = -(r-1)z + r (r-1) = -(r-1)(z-r)$, and so  $rf$ is decreasing as a function of $x$ (keeping $z $ fixed, but at least $r+1$), and thus it suffices to evaluate $rf$ at at $z = (10r- 8)/2 +1 - r = 4r-3$.
$$\eqalign{
rf(4r-3) \left|_{x=10r-8}\right. & = r(4r-3)^2 -(r-1)(9r-8) \(4r -3\) + r(r-1) (10r-8)\cr  
&= -10r^3 +  53 r^2 - 66r + 24\cr
 }$$
This is obviously  negative for $r \geq 6$, but by direct computation, we also see it is negative if $r = 4$ or $5$. 

Next,  $4rf(2r+1) = - 8r(r-1)x + 4r (6r^2 + 3r + 2)$. This is negative if $x > (6r^2 + 3r + 2)/2(r-1) = 3r + 3 + 4/(r-1)$, and this is true if $x \geq 5r$ and $ r > 2$. \qed

\noindent {\it Proof of theorem \Aone, completed.}   We may assume $(x,r) = 1$. By Lemma \Asix(d) and (c), we need only show $S(x,r) \cap [r+1, x-r] = \emptyset$. By \Asix(e) and Lemma \Asev, it is sufficient (assuming $ x \geq r^2 \vee 10$) that $S(x,r) \cap [2r+1, x/2 + 1 - r] = \emptyset$.

The set of $z \in \R$ \st $f(z) \leq 0$ (with $x$ and $r$ fixed) is (possibly empty) closed interval. Thus by Lemmas \Aele\ and \Aeig, if $x \geq 10r - 8$ and $r\neq 2,3$, then $S(x,r) \cap [2r+1, x/2 + 1 -r ] = \emptyset$, and we are done in this case. 

Next, we observe that $x/2 + 1-r \leq 2r$ if and only if $x \leq 6r-2$. By Corollary \Asev,   the Theorem follows here. 

Hence we are reduced to the following situations, with $(x,r) = 1$ and $x \geq r^2+1$:
\item{(i)} $r = 2$;
\item{(ii)} $r = 3$;
\item{(iii)} $r \geq 4$  and $6r-1 \leq x \leq 10r-8$.
 \vskip 2pt
\noindent If $r \geq 9$, then $r^2 + 1 > 10r-8$, so we are done. 

\vskip5pt
\noindent (i) {\it When $r = 2$,} we have that $x$ is odd, and by  direct calculation, we can assume that $x \geq 21$. By Lemma \Asix(f), it suffices to show $S(x,2) \cap [2r+1, (x-9)/2] = \emptyset$; for this, it suffices to show $f$  (or $rf$) is negative at the two endpoints; $f(5) = 25 - 5 (x-2) + 2x = -3x + 35$, so $x \geq 13$ is sufficient, and $4f((x-9)/2) = - 3x + 63$, which is non-positive, as $x \geq 23$. 

\vskip 5pt 
\noindent (ii) {\it Now suppose that $r = 3$,}   and $(x,3) = 1$; by calculation, we may assume $x \geq 31$. Then $rf(x/2 + 1- r) = 3f((x-5)/2)  $, and this evaluates to $(-x^2 + 22 x + 15)/4$); this is negative, as $x \geq 23$. On the other hand, $f(2r+1) = f(7) = 49 + 14 - 8x/3$, and this is non-positive if $x \geq 189/8 = 23\slfrac{5}8$, that is, $x \geq 24$ (since $(x,3) =1$, $x \geq 25$). We are left with just the case that $x = 23$. Here $x/2 + 1 -r =  19/2$; hence we need only check that $\brcs{7,8, 9} \cap S(23,3) = \emptyset$, which is easy: for $z = 7$, set $t = 3$; for $z = 8$, set $t = 5$; and for $z = 9$, set $t = 5$. 
\vskip 2.5 pt

 \noindent (iii) {\it Subcase (a)  $6r+5 \leq x \leq	8r-2$ for $r \geq 5$, and $6r+ 6 \leq x \leq 8r-2$ if $r=4$.} We have $f(2r+1) \leq 0$ if $x \geq 5r$, and a direct computation yields $f(2r+3) \leq 0 $ if $x \geq 6r+1 + (22r+ 3)/(r^2 +2r -3):= x_0 (r)$. The values of $\ceil{x_0 (r)}$ for $r = 4, \dots, 8$ are respectively, $\(30,35,40,46,52\) = \brcs{6r + 6, 6r+5, 6r+5,6r + 5, 6r + 5, 6r+4 }$.  Thus if $6r+5 \leq x  $ for $r \geq 5$ or $x \geq 30$ if $r = r$, we have that $f$ is nonnegative on the interval $[2r+1, 2r+3]$, and so $[2r+1, 2r+3] \cap S_0(x,r) = \emptyset$. By \Asix(d), we have $[2,\dots, 3r] \cap S_0 (x,r) = \emptyset$ (this requires $r \geq 4$). 
 We will be done with this case if $3r \geq x/2 + 1 - r$, that is, $x \leq 8r -2$, and thus with the situation that $6r + 5 \leq x \leq 8r-2$ and $r \geq 5$ or $x=30$ and $r = 4$. 
 
\vskip 2.5 pt
  {\it Subcase (b) $6r-1 \leq x \leq 6r+ 4$ for $r \geq 5$ and $23 \leq x \leq 29$ for $r = 4$.} We will show that $[2r+1, 2r+3]\cap S_0(x,r) = \emptyset$ again in this case, yielding validity of the theorem for $6r-1 \leq x \leq 8r-2$. This is by direct computation.

\vskip 1.5 pt
 \noindent{$(r=4)$.} For $z = 9,10,11$ and $x = 23,25, 26, 27, 29$ ($24$ and $28$ are excluded as $r|x$, and this was already done), the corresponding choices for $t $ will work: for $z =9$, they are $t=2,5,5,5,3,$ respectively; for  $z = 10$, we take $t = 2,4,5,5,5$; and for $z=11$, we choose $t = 2,4,5,5,5$ again.

\vskip 1.5 pt
 \noindent{$(r=5)$.} For $z = 11,12,13$ and $x = 29,31,32,33,34$ ($30$ is excluded as $r|x$e), the corresponding choices for $t $ will work: for $z =11$, set  $t=5$ (all for each of the five values of $x$);    $z = 12$, we take $t = 4,5,5,5,5$; and for $z=11$, we choose $t = 2,2,2,2,8$. 

\vskip 1.5 pt

 \noindent{$(r=6)$.} For $z = 13,14,15$ and $x = 37, 38, 39, 40$ (we are assuming $x \geq r^2$, and $r$ does not divide $x$), the corresponding choices for $t $ will work: for each of the three values of $z$ and the four values of $x$, except $(15,37)$, we take $t = 5$; if $ z = 15$ and $x = 37$, set $t = 4$. 

\vskip 1.5pt

 For $r \geq 7$, $6r +4$ falls below the threshhold $x \geq r^2 +1$ and there is nothing to do. 
 
 \vskip 1.5pt Subcase (c) $8r-1 \leq x \leq 10r-9$.  Assume this.  From (b) above, as $8r-1 \geq 6r+5$, we have that $\brcs{2r+1, 2r+2,2r+3} \cap S_0 (x,r) = \emptyset$, hence $\brcs{2,\dots, 3r} \cap S_0(x,r) = \emptyset$.  Now we show that  $\brcs{3r+1, 3r+2, 3r+3}$. 
 
 If $r = 4$, then $x = 31$, and $z = 13, 14, 15$, then $t = 2$ will work. Now $15 \geq x/2 + 1 - 4$ is implied by  $x \leq 34$,  and we are done.  
 
 If $r \geq 5$, we show that $f(3r + a) \leq 0$ when $a = 1$ and $a = 4$; this will entail $f$  is nonpositive on $[3r +1, 3r+4]$, and thus $[3r+1, 3r+4]\cap S_0(x,r) = \emptyset$. This will imply $[2,\dots , 4r] \cap S_0(x,r) = \emptyset$, and sufficient for $4r \geq x/2 +1-r$ is $x \leq 10r-2$, and the latter exceeds $10r-9$, completing the proof. 
 
 An elementary computation yields that $f(3r + a) \leq 0$ iff 
$$
 x \geq  6r + \frac{a+9}2 + \frac{\(\frac{a^2}2 + \frac 72 a + 9\)r - a(a+9)}{2r^2 + (a-2)r -a}.
$$
For $a = 1$, the respective values of the ceiling of the right side for $r = 5,6,7,8$ are $37, 43,48,54$, each of which is less than $8r-1$. For $a =4$, the respective values  are $39, 45, 51, 57$, each of which is bounded above by $8r-1$. Hence $f$ is nonpositive on $[3r+1,3r+4]$, and we are done. \qed